\documentclass[11pt]{amsart}
\usepackage{amssymb,amstext}
\newtheorem{theorem}{Theorem}[section]
\newtheorem{lemma}[theorem]{Lemma}

\theoremstyle{definition}
\newtheorem{definition}[theorem]{Definition}
\newtheorem{prop}[theorem]{Proposition}
\newtheorem{diag}[theorem]{Diagram}
\newtheorem{cor}[theorem]{Corollary}

\newtheorem{example}[theorem]{Example}

\theoremstyle{remark}
\newtheorem{remark}[theorem]{Remark}

\numberwithin{equation}{section}

  \newcommand{\Q}{\mathbb Q}
  \newcommand{\R}{\mathbb R}  
  \newcommand{\N}{\mathbb N}
  \newcommand{\C}{\mathbb C}
 \newcommand{\Z}{\mathbb{Z}}
 \newcommand{\F}{\mathbb{F}}
 \newcommand{\bA}{\mathbb{A}}
 \newcommand{\bT}{\mathbb{T}}
 \newcommand{\oA}{\overline{\bA}}

 \newcommand{\rk}{\mbox{\rm rk}\,}
 
 \newcommand{\T}{\mathcal T}

 \newcommand{\D}{\mathcal D}
\newcommand{\im}{\mbox{\rm im\,}}
\newcommand{\coker}{\mbox{\rm coker\,}}
\newcommand{\rank}{\mbox{\rm rank\,}}
\newcommand{\uA}{\underline{A}}
\newcommand{\uT}{\underline{T}}
\newcommand{\uD}{\underline{D}}
\newcommand{\uM}{\underline{M}}
\newcommand{\uN}{\underline{N}}
\newcommand{\uP}{\underline{P}}
\newcommand{\uX}{\underline{X}}
\newcommand{\uY}{\underline{Y}}
\newcommand{\uC}{\underline{C}}

\newcommand{\bra}[1]{\stackrel{#1}{\longrightarrow}}

\newcommand{\W}{{\mathcal W}}
\newcommand{\JJ}{\mathcal P}
\newcommand{\RR}{\mathcal R}
\newcommand{\erz}[1]{\langle #1 \rangle}
\newcommand{\dpe}{{\dim \! V}}

\newcommand{\Gammaup}{\Gamma}
\newcommand{\Gammaperp}{\Gamma^\perp}
\newcommand{\Gammapar}{\Gamma^\|}
\newcommand{\piperp}{\pi^\perp}
\newcommand{\pipar}{\pi^\|}
\newcommand{\Gup}{\Gammaup}

\newcommand{\Eup}{E}
\newcommand{\Epar}{E^\|}
\newcommand{\Eperp}{E^\perp}
\newcommand{\Wup}{D}
\newcommand{\Torus}{\mathbb T}
\newcommand{\MP}{\Omega}
\newcommand{\co}{\colon\, }
\newcommand{\can}{almost canonical} 
\newcommand{\met}{pattern metric}
\newcommand{\onto}{\longrightarrow \!\!\!\!\!\!\! \rightarrow}

\begin{document}
\title{Integral cohomology of rational
  projection method patterns}

\author{Franz G\"ahler}
\address{Faculty of Mathematics, University of Bielefeld, D-33615 Bielefeld, Germany.}
\email{gaehler@math.uni-bielefeld.de}

\author{John Hunton}
\address{The Department of Mathematics, University of
Leicester, Leicester, LE1 7RH, England.}

\email{j.hunton@mcs.le.ac.uk}

\author{Johannes Kellendonk}
\address{Universit\'e de Lyon, Universit\'e Claude Bernard Lyon 1,
Institut Camille Jordan, CNRS UMR 5208, 43 boulevard du 11 novembre
1918, F-69622 Villeurbanne cedex, France.}

\email{kellendonk@math.univ-lyon1.fr}

\date{\today}

\begin{abstract}
We study the cohomology and hence $K$-theory of the aperiodic tilings formed by the so called \lq cut and project\rq\ method, i.e., patterns in $d$ dimensional Euclidean space which arise as sections of higher dimensional, periodic structures. They form one of the key families of patterns used in quasicrystal physics, where their topological invariants carry quantum mechanical information. Our work develops both a theoretical framework and a practical toolkit for the discussion and calculation of their integral cohomology, and extends previous work that only successfully addressed rational cohomological invariants. Our framework unifies the several previous methods used to study the cohomology of these patterns. We discuss explicit calculations for the main examples of icosahedral patterns in $\R^3$ -- the Danzer tiling, the Ammann-Kramer tiling and the Canonical and Dual Canonical $D_6$ tilings, including complete computations for the first of these, as well as results for many of the better known 2 dimensional examples.
\end{abstract}

\maketitle

\bibliographystyle{amsalpha}

\tableofcontents
\section{Introduction}
This work considers one of the key families of aperiodic patterns used in quasicrystal physics. We develop both a theoretical framework and a practical toolkit for the discussion and calculation of the integral cohomology and $K$-theory of these patterns. Our work extends previous results which successfully addressed only their {\em rational\/} cohomology \cite{FHKcmp, FHKmem, Kal} and it provides a unified treatment of the two apparently distinct approaches \cite{FHKcmp, FHKmem} and \cite{Kal} studied so far in the literature. The patterns we consider are point patterns in some $d$-dimensional Euclidean space $\R^d$ that arise as sections of higher-dimensional, periodic structures, variously known as  {\em model sets}, {\em cut \& project\/} patterns or just {\em projection patterns\/} \cite{moody}. By a standard equivalence, such point patterns may also be considered as tilings, coverings of $\R^d$ by compact polyhedral sets meeting only face to face. The Penrose tiling in 2 dimensions is perhaps the best  known example, but the class is huge (indeed, it is infinite) and today forms the principal set of geometric models for physical quasicrystals; see, for example, \cite{SD}

To any point pattern or tiling $P$ in $\R^d$ a topological space associated to $P$, called the {\em hull\/} or {\em tiling space\/} $\MP$ of $P$, may be constructed. In short, this is a moduli space of patterns locally equivalent to $P$. Under standard assumptions (certainly satisfied by the class of patterns we consider), $\MP$ is a compact, metrisable space, fibering over a $d$-torus with fibre a Cantor set \cite{FHKmem, SW}. Much progress during the last 20 years or more in the study of aperiodic patterns has developed through the study of these spaces, which can by analyzed via standard topological machinery such as cohomology or $K$-theory. Major results include the Gap Labeling Theorem \cite{Bellissard, BBG, BHZ, BKL, BO, KP}, the deformation theory of tilings \cite{CS, Keldef, SW} and the work on exact regularity of patterns and the homological Pisot conjecture \cite{BBJS, SadER}. For a short introduction to the topology of tiling spaces and some of the geometric and physical benefits of understanding their cohomology, we direct the reader also to \cite{Sadbook}.

It is a general truth that by writing any tiling space $\MP$ as a Cantor bundle over a $d$-torus, one can realise the Cech cohomology of $\MP$ as the group cohomology of $\Z^d$ with coefficients derived from the structure of the fibre and the holonomy of the bundle. In general, however, one has little hold over either the fibre or the holonomy, but, as was realised in \cite{FHKcmp, FHKmem}, there is a large class of projection tilings for which a practical description can be obtained. This class contains the so-called canonical projection tilings, and was later called the class of {\em almost canonical\/} tilings \cite{AJ}; we present them formally in the next section. (It is interesting to note that this is also the class of tilings whose asymptotic combinatorial complexity can be easily obtained \cite{AJ}.)

In \cite{FHKcmp, FHKmem} Forrest, Hunton and Kellendonk effectively provided a method for the computation of the rational cohomology of the spaces $\MP$ of almost canonical tilings. A related, but non-commutative approach, describing the $K$-theory of crossed product algebras associated to these tilings, was given by Putnam in \cite{Putnam}. Results similar to \cite{FHKmem} for a smaller class of projection tilings, produced from an apparently rather different perspective, were obtained by Kalugin in \cite{Kal} who gave a shape equivalent approximation to $\MP$ by a finite CW complex (though that terminology was not used in \cite{Kal}). 

However, a key feature of the interpretations of all these works at the time was the assumption that the cohomology and $K$-theory of these pattern spaces would be free of torsion, and thus integral computation would follow from just working with rational coefficients and counting ranks of vector spaces. This turned out not to be the case (and, unfortunately, some statements about and referring to the torsion freeness of cohomological or $K$-theoretic invariants in \cite{FoHu, FHKcmp, FHKmem} are wrong). This was shown, for example, by G\"ahler's counterexamples \cite{FGBanff} obtained through extensive machine computation for certain 2 dimensional patterns which arise as both projection and substitution tilings. The substitution structure allowed a yet further approach to computation via the method of Anderson and Putnam \cite{AP}, though even for relatively modest 2 dimensional examples this method is stretched to the limit of accessible computation. Nevertheless, examples computed, in particular the T\"ubingen Triangle Tiling (TTT) \cite{BKSZ,KSB}, demonstrated that the integral cohomology could be far more complicated than had previously been thought, and this formed the stimulus of our work here. We understand that the existence of torsion, not appreciated at the
time when \cite{BBG, BKL, BO, KP} were written, may cause problems with some of the arguments used in the published proofs of the Gap Labeling Theorem.

Given the consequent complexity of the cohomology $H^*(\MP)$, its  complete description for projection method patterns is beyond the scope of the techniques of any of \cite{AP, FHKcmp,FHKmem, Kal} for all but the simplest examples. 

In this paper we present techniques to address this. In Section \ref{HomAlg} we introduce a set of ideas from homological algebra that can be applied for discussing the bundle structures associated to these patterns. As a further  consequence, the generality of the framework developed allows us to unify the approaches of both \cite{FHKcmp,FHKmem} and \cite{Kal}, and this point  has computational advantages when we turn in the final section to the discussion of the more complex examples.

In Section \ref{GeoReal} we give a geometric interpretation of almost canonical projection patterns whose cohomology is finitely generated and which satisfy one further assumption. This is inspired by and is an analogue of a certain key assumption made in Kalugin's approach \cite{Kal}. Patterns which enjoy this geometric interpretation we term {\em rational projection method patterns}; they  form the central class for which we compute integer cohomology in the final section. Section \ref{GeoReal} ends with a complete description of the cohomology of a rational projection pattern in terms of data coded in the cohomology of an inclusion of a certain finite CW complex $\bA$ in an ambient torus $\bT$.

It is these two new ingredients, the geometric interpretation of Section \ref{GeoReal} and the homological framework of Section \ref{HomAlg}, which give us tools to analysis integer cohomologyy for examples beyond the ready scope of any of the previous works in the field.

The final sections of the paper turn to the actual computation of examples. The complexity of the computation of the cohomology of a projection pattern increases with the so-called {\em codimension\/} of the pattern. In Section \ref{Sect5} we give a complete description of the cohomology of rational projection patterns of codimension 1 and 2, together with details of many of the main examples and an outline of the machine methods used to compute them. Strictly speaking, the results of this section are accessible with the older techniques of  \cite{FHKcmp,FHKmem, Kal}, but the section provides the necessary foundation for the new and more complex work of Section \ref{Sect6} which considers the codimension 3 examples and, briefly, the cohomology and $K$-theory of general codimension rational patterns. We note that the physically interesting rational projection patterns (i.e., those in dimension up to 3) arise only from codimension 1, 2 or 3 schemes. We compute explicitly the cohomology of the Danzer tiling \cite{Dan89}, and much of the cohomology of three other  3 dimensional, icosahedral patterns, those of Ammann-Kramer  \cite{AmmannKramer}, the canonical $D_6$ and dual canonical $D_6$ patterns \cite{KP95}.

Some of these results and ideas were announced in \cite{GHK} (though the reader should note that there are some errors in the computation of the torsion component of $H^3(\Omega)$ published in \cite{GHK} -- see Section \ref{sec:ex3d} for details), but the framework and techniques presented here have developed considerably since that note.

\smallskip\noindent{\bf Acknowledgements.} The first author was supported by the German Research Council (DFG) within
the CRC 701, project B2.  The second author acknowledges the support of study leave granted by the University of Leicester, and the hospitality of Universit\'e de Lyon. The third author acknowledges the financial support of the ANR {\em SubTile}.

\section{Projection patterns, their spaces and cohomology}\label{Sect2}
We begin by describing the types of patterns we consider, and in so doing set up our notation. The contents of this section are mostly a brief summary of the set-up and foundational results of \cite{FHKcmp, FHKmem}; the reader should consult those sources for further detail and discussion. We start by listing the data needed for a model set, or cut and project pattern.

\begin{definition}\label{CandP} A {\em cut \& project scheme} consists of a euclidean space $\Eup$ of dimension $N$ containing a discrete cocompact abelian group (or {\em lattice}) $\Gammaup$. There is a direct sum decomposition  $\Eup=E^\|\oplus E^\perp$ with associated projections $\pi^\|\co  \Eup\to E^\|$ and $\pi^\perp\co  \Eup\to E^\perp$. We assume $E^\|$ and $E^\perp$ are in {\em total irrational position\/} meaning that $\pi^\|$ and $\pi^\perp$ are one to one and with dense image of the lattice $\Gammaup$. Denote by $d$, respectively $n$, the dimensions of $E^\|$ and $E^\perp$, so $N=d+n$. We call $d$ the {\em dimension\/} of the scheme, and $n$ its {\em codimension}. Finally, we have also an {\em acceptance window\/} or {\em atomic surface\/} $K$, a finite union of compact non-degenerate polyhedra in $E^\perp$. We denote by $\partial K$ the boundary of $K$, which consists of a finite union of $(n-1)$-dimensional faces $\{f_i\}$.
\end{definition}

For convenience we denote by $\Gammapar$ and $\Gammaperp$ the images $\pi^\|(\Gammaup)$ and $\pi^\perp(\Gammaup)$. These are both rank $N$ free abelian subgroups of $E^\|$ and $E^\perp$ respectively. 

\begin{definition} Given a cut \& project scheme, we define the associated point pattern $P$ as the set of points in $\Epar$ 
$$P=\{\pi^\|(\gamma)|\gamma\in\Gammaup:\pi^\perp(\gamma)\in K\}$$
or equivalently  as
$$P=\Epar\cap(\Gammaup-K).$$
\end{definition}

There are a number of variations in the way cut and project patterns can be viewed. In \cite{FHKmem}  the viewpoint was taken that these patterns arise as projections of point patterns within strips $E^\|+K$.  
Kalugin in \cite{Kal} uses the section method by means of which these patterns arise as intersections between $\Epar$ and a $\Gammaup$-periodic arrangement of sets. In \cite{FHKcmp} the dual method using Laguerre complexes was adopted, which is more elegant for some tilings such as the Penrose tilings. The reader can consult Moody's work, for example \cite{moody}, for a wide ranging discussion of these patterns.

A cut and project scheme in fact defines a whole parameterised family of point patterns in $E^\|$.

\begin{definition}
For each point $x\in \Eup$ define the point set
$$\begin{array}{rcl}P_x &=&\{\pi^\|(\gamma)|\gamma\in\Gammaup:\pi^\perp(\gamma+x)\in K\}\\
&=&  \Epar\cap(\Gammaup+x-K).
\end{array}$$
\end{definition}
Note that the pattern $P_x$ depends only on the class of $x$ in $\Eup/\Gammaup=\Torus$, an $N$-torus. In fact $P_x=P_y$ if and only if $x-y\in\Gup$.

\begin{definition} We define the set $S$ of {\em singular points\/} in $\Eup$ by
$$ S =\{ x\in\Eup : \pi^\perp(x)\in \partial K+\Gammaperp\} = E^\| + \Gammaup + \partial K\,.$$
Denote by $N\!S$ its complement, the set of {\em nonsingular points}.
\end{definition} 

It is well known that, for any $x$, the pattern $P_x$ is {\em aperiodic}, i.e., that $P_x=P_x+v$ only if $v$ is the zero vector, and is of {\em finite local complexity}, meaning that, up to translation, for each $r>0$ there are only a finite number of local configurations of radius $r$ in $P_x$. If $x\in N\!S$ then $P_x$ satisfies the additional property that for each finite radius $r$ there is a number $R$ such that any radius $r$ patch of $P_x$ occurs within distance $R$ of any given point of $\Epar$, a property known as {\em repetitivity}. We note further that  if $x$ and $y$ are both nonsingular points, then the patterns $P_x$ and $P_y$ are locally indistinguishable in the sense that each compact patch of one pattern occurs after translation as a patch in the other. Although these are important properties and motivate interest in understanding and characterising cut and project patterns, they will not generally play a very explicit role in the work which follows, though they implicitly account for many of the topological properties of the space $\MP$ we will shortly introduce and is the main topic of the article. Again, see \cite{moody} for further introduction and discussion of these properties.

The cohomology of point patterns which we investigate here is the Cech cohomology of an associated pattern space. Suppose for simplicity that $0\notin S$.

\begin{definition}
The pattern space $\MP$ of $P=P_0$ is the completion of the translates of $P$
with respect to the {\em \met\/}, defined on two subsets $P,Q\subset E^\|$ by
$$ d(P,Q) = \inf\left\{\frac{1}{r+1}\left|
\begin{array}{rl}
&\!\!\!\!\mbox{ there exists }x,y\in B_\frac{1}{r} \mbox{ with}\\
&\big(B_r\cap (P-x)\big)\cup \partial B_r = \big(B_r\cap (Q-y)\big)\cup \partial B_r 
\end{array}\right. \right\}.$$ 
Here $B_r$ is the closed ball around $0$ of radius $r$ in $E^\|$. In essence this metric is declaring two patterns to be close if, up to a small translation, they are identical up to a long distance from the origin. The precise values of this metric will not be important in what follows, but rather the topology it generates.
\end{definition}

It can readily be shown that the space $\MP$ contains precisely those point patterns which are locally indistinguishable from $P$. As $P_x=P_y$ if and only if $x-y\in\Gup$, $\MP$ can also be seen as the completion of $q(N\!S)\subset\Torus$ with respect to the \met, where $q\co \Eup\to\Torus$ is the quotient $\Eup\to\Eup/\Gammaup$. Furthermore, the same space $\MP$ is obtained on replacement of $P$ in the previous definition by $P_x$ for any nonsingular $x$.

\begin{definition}
The cohomology of a projection method pattern $P$ is the Cech cohomology of the associated space $\MP$. We shall denote this $H^*(\MP)$ when we are considering coefficients in $\Z$, and by $H^*(\MP;R)$ when we take coefficients in some other commutative ring $R$.
\end{definition}

Note that the \met\ is not continuous in the euclidean topology of the parameter space $q(N\!S)\subset\Torus$ but conversely, the euclidean metric on $\Torus$ is continuous with respect to the \met. Therefore there is a continuous map 
$$\mu\co  \MP\to\Torus,$$ 
in fact a surjection, such that each non-singular point has a unique pre-image. Since $q(N\!S)$ is large in a topological sense (it is a dense $G_\delta$-set) and in the measure sense (it has full Lebesgue measure) $\mu$ is called almost one to one. See \cite{FHKmem} for a full discussion.

\begin{definition}  We shall call the cut \& project scheme (and its corresponding patterns) {\em \can\/} if 
for each face $f_i$ of the acceptance domain, the set $f_i+\Gamma^\perp$ contains the affine space spanned by $f_i$.
\end{definition}

We assume throughout this paper that our scheme and patterns are \can. From the constructions of \cite{FHKmem} it can be shown \cite{irving} that for the patterns of Definition \ref{CandP} this is a necessary (but certainly not sufficient) condition for the Cech cohomology $H^*(\MP)$ to be finitely generated.

This definition is equivalent to saying that there is a finite family of $n-1$ dimensional affine subspaces
$$\W=\{W_\alpha\subset E^\perp\}_{\alpha\in I_{n-1}}$$ such that 
$$ S = E^\| + \Gammaperp + \bigcup_{\alpha\in I_{n-1}} W_\alpha .$$
Note that we have some freedom to choose the spaces $W_\alpha$: replacing  
$W_\alpha$ by $W_\alpha-\gamma$ for some $\gamma\in\Gammaperp$ does not change the singular set $S$. We will always assume that $\W$ has the least number of elements possible, which means that from every $\Gammaperp$-orbit we have only one representative. 

\begin{definition}
Suppose the cut \& project scheme is \can, and we have chosen some such family of subspaces $\W$. Call an affine subspace $W_\alpha+\gamma\subset E^\perp$, for any $\alpha\in I_{n-1}$ and $\gamma\in\Gammaperp$ a {\em singular space}. Clearly the set of all singular spaces is independent of the particular finite family $\W$ chosen. 
\end{definition}

Intersections of $\Gammaperp$-translates of singular spaces may be empty, but if not they yield affine subspaces of lower dimension. We shall call all affine spaces arising in this way singular spaces as well.  Note that $\Gamma$ acts on the set of all singular spaces by translation; if $\gamma\in\Gamma$ and $W$ is a singular space of dimension $r$, then so is $\gamma\cdot W=W+\pi^\perp(\gamma)$. The {\em stabilizer\/} $\Gamma^W$ of a singular space $W$ is defined as the subgroup of $\Gamma$ given by
$\{\gamma\in\Gammaup | W-\pi^\perp(\gamma) = W\}$. Note that the stabilizers of singular spaces which differ by a translation coincide.

The cohomology groups $H^*(\MP)$ depend on the geometry
and combinatorics of the intersections of the singular spaces and the action of $\Gamma$ on them. It will therefore be useful to develop notation for these concepts. Recall that $I_{n-1}$ indexes the set of orbit classes of all $(n-1)$-dimensional singular spaces.

\begin{definition}
\begin{enumerate}
\item For each $0\leqslant r<n$, let $\JJ_r$ be the set of all singular $r$-spaces. Denote the orbit
space under the action by translation $I_r=\JJ_r/\Gamma$.
\item The stabilizer $\Gamma^W$ of a singular
$r$-space $W$ depends only on the orbit class $\Theta\in
I_r$ of $W$ and we will also denote it $\Gamma^\Theta$. 
\item Suppose $r<k<n$ and pick some $W\in\JJ_k$ of orbit class $\Theta\in I_k$.
Let $\JJ_r^{W}$ denote $\{U\in\JJ_r|U\subset W\}$, the set of singular $r$ spaces lying in $W$. 
Then $\Gamma^\Theta$ acts on $\JJ_r^{W}$ and we write 
$I_r^{\Theta}=\JJ_r^{W}/\Gamma^\Theta$, a set which depends only on the class $\Theta$ of $W$. Thus
$I_r^{\Theta}\subset I_r$ consists of those orbits of singular
$r$-spaces which have a representative that lies in a singular
$k$-space of class $\Theta$.
\item  Finally we denote the cardinalities of these sets by $ L_r =
|I_r|$ and $L^\Theta_r = |I^\Theta_r|$.
\end{enumerate}
\end{definition}

We recall some of the main results of \cite{FHKmem}. 

\begin{theorem}\label{L0stuff}
\begin{enumerate}
\item$L_0$ is finite if and only if $H^*(\MP)$ is finitely generated as a graded abelian group.  {\em \cite{FHKmem}, Theorems IV.2.9 \& V.2.5.} 
\item If $L_0$ is finite then all the $L_r$ and $L_r^\Theta$ are
  finite as well, and $\nu=N/n$ is an integer. Moreover, $\mbox{\rm rank}\,\Gamma^U=\nu\cdot\dim(U)$  for any singular space $U$ if and only if $L_0$ is finite.  {\em \cite{FHKmem}, Lemma V.2.3 \& Theorem IV.6.7, and \cite{AJ}.} 
\end{enumerate}
\end{theorem}

\section{Homological algebra for cut and project schemes}\label{HomAlg}
\subsection{$\mathcal C$-topes and complexes}

We assume we have a \can\ cut \& project scheme, with associated $(n-1)$-dimensional singular spaces $\JJ_{n-1}=\{W_\alpha+\Gammaperp\}_{\alpha\in I_{n-1}}$ in $E^\perp$. The geometry and combinatorics of these spaces give rise to a $\Gamma$-module $C_n$ key to our work on $H^*(\MP)$. The module $C_n$, and associated objects given by the lower dimensional singular spaces, were first introduced in \cite{FHKcmp, FHKmem} where the equivalences
\begin{equation}\label{FundEquiv}
H^s(\MP)\cong H^s(\Gamma;C_n)\cong H_{d-s}(\Gamma;C_n)\,.
\end{equation}
were shown (eg, \cite{FHKcmp} corollaries 41, 43). Here the latter two groups are the group cohomology, respectively group homology, of $\Gamma$ with coefficients in $C_n$.

We outline the proof of these equivalences in the Appendix, and complete details can be found in \cite{FHKmem}, but for now we recall the definition of $C_n$ and associated modules, and develop further related algebraic tools.

\begin{definition}\label{C-tope}
Call a  $\mathcal C$-tope any
compact polyhedron $J$ in $\Eperp$ whose boundary belongs to some union 
$\bigcup_{W\in A}W$, where $A$ 
is a finite subset of $\JJ_{n-1}$. As on singular spaces, 
$\Gamma$ acts on the set of $\mathcal C$-topes by translation, $\gamma\cdot J = J-\piperp(\gamma)$. 
Each connected component of the window $K$ is a $\mathcal C$-tope and, in fact, all $\mathcal C$-topes occur as components of finite unions of finite intersections of $\Gammaperp$-translates of $K$.  
\end{definition}

Let $C_n$  be the $\Z\Gamma$-module generated by indicator functions on $\mathcal C$-topes, and for $r<n$ let $C_r$ be the $\Z\Gamma$-module generated by
indicator functions on $r$-dimensional facets of $\mathcal C$-topes. 
In particular, $C_n$ can be identified with $C_c(\Eperp_c,\Z)$, the $\Z$-module of compactly supported $\Z$-valued functions on $\Eperp$ with discontinuities only at points of $\JJ_{n-1}$. 

The set $\JJ_{n-1}$ of all singular $(n-1)$ spaces is dense in $\Eperp$. It will be useful to view $\JJ_{n-1}=\cup_i\JJ_{n-1}(i)$ where $\JJ_{n-1}(1)\subset \JJ_{n-1}(2)\subset\cdots\subset \JJ_{n-1}(i)\subset\cdots$ is an increasing sequence of locally finite collections of singular $(n-1)$-spaces. Write also $\JJ_r(i)$ for the singular $r$-spaces occurring as intersections of the elements of $\JJ_{n-1}(i)$. For $r<n$ denote by $C_r(i)$ the $\Z$-module of compactly supported $\Z$-valued functions on the singular $r$-spaces in $\JJ_r(i)$ with discontinuities only at points of $\JJ_{r-1}(i)$, and for $r=n$ write  $C_{n}(i)$ for the $\Z$-module of compactly supported $\Z$-valued functions on $\Eperp$ with discontinuities only at points of $\JJ_{n-1}(i)$. Clearly there are inclusions $C_r(i)\to C_r(i+1)$ and this construction yields 

\begin{lemma}\label{limit}$$C_r=\displaystyle{\lim}_{i\to\infty} C_r(i)\,.$$\qed\end{lemma}

These modules form a complex of $\Z\Gamma$-modules with $\Gamma$-equivariant boundary maps
\begin{equation}\label{complex}
0 \to C_n\stackrel{\delta}{\rightarrow} C_{n-1}\stackrel{\delta}{\rightarrow} \cdots \stackrel{\delta}{\rightarrow} C_0 \stackrel{\epsilon}\rightarrow \Z \to 0,
\end{equation}
$\delta$ being induced by the cellular boundary map on $\mathcal C$-topes and $\epsilon$ the augmentation map defined as follows. The module $C_0$ is generated by indicator functions on 0-dimensional singular spaces; denote such a function by $1_p$ for some $p\in\JJ_0$. Then $\epsilon$ is given by $\epsilon(1_p)=1$.

\begin{lemma}(\cite{FHKcmp} Prop 61)
The sequence of $\Z\Gamma$ modules (\ref{complex}) is exact.
\end{lemma}

\smallskip\noindent{\bf Sketch proof.} First note that the corresponding sequence
$$0 \to C_n(i)\stackrel{\delta}{\rightarrow} C_{n-1}(i)\stackrel{\delta}{\rightarrow} \cdots \stackrel{\delta}{\rightarrow} C_0(i) \stackrel{\epsilon}\rightarrow \Z \to 0$$
is the augmented cellular chain complex of the space $\Eperp$ with cellular decomposition given by the family of hyperplanes $\JJ_{n-1}(i)$. It is exact since $\Eperp$ is contractible. The result follows by taking the direct limit as $i\to\infty$: exactness is preserved by direct limits.\qed

\smallskip
It will be useful to have a homological interpretation of the modules $C_r$ and this will follow from the cellular structures induced by the $\JJ_{n-1}(i)$ as in the proof of the last lemma. For convenience we shall denote also by $\JJ_r(i)$, etc, the subspace of $\Eperp$ consisting of the union of the affine subspaces in this set.

\begin{lemma}\label{CWcomp}
$$\begin{array}{rl}
\displaystyle{\lim}_{i\to\infty}  H_m(\Eperp, \JJ_{n-1}(i))&=\left\{
\begin{array}{ll}C_n&\mbox{if m=n,}\\
0&\mbox{otherwise;}\\
\end{array}\right.\\
\displaystyle{\lim}_{i\to\infty}  H_m(\JJ_r(i),\JJ_{r-1}(i))&=\left\{
\begin{array}{ll}C_r&\mbox{if m=r,}\\
0&\mbox{otherwise.}\\
\end{array}\right.
\end{array}$$
Here $H_*(X,Y)$ denotes the relative homology of the pair $Y\subset X$.
\end{lemma}

\smallskip\noindent{\bf Proof.} If $X$ is a CW complex with $r$-skeleton $X^r$ (i.e., the union of all cells of dimension at most $r$), then $X^r/X^{r-1}$ is a one point union of $r$-spheres, in one-to-one correspondence with the $r$-cells of $X$. Thus $H_r(X^r,X^{r-1})=H_r(X^r/X^{r-1})$ is the $r^{\rm th}$ cellular chain group for $X$ while $H_t(X^r,X^{r-1})=0$ for $t\not=r$. As $\JJ_r(i)$ is the $r$-skeleton of $\Eperp$ with CW structure given by the $\JJ_{n-1}(i)$, the lemma follows by taking limits as $i\to\infty$.\qed

\smallskip Finally we note the following decomposition results for the lower $C_r$. Full details can be found in \cite{FHKmem} Lemma V.3.3 and Corollaries V.4.2, \& V.4.3. For $r<n$ and $\alpha\in I_r$, if $W$ is a singular $r$-space representative of the orbit  indexed by $\alpha$, write $C_r^\alpha$ for the $\Z[\Gamma^\alpha]$-module of $\Z$-valued functions on $W$ with discontinuities where $W$ meets transversely the singular spaces $\JJ_{n-1}$. Similarly, for $r<k<n$ if $V$ is a singular $k$-space of orbit  class $\alpha\in I_k$, and $W$ is a singular $r$-space in $V$ of orbit class $\psi\in I_r^\alpha$, write $C_r^{\alpha,\psi}$ for the $\Z[\Gamma^\psi]$-module of $\Z$-valued functions on $W$ with discontinuities where $W$ meets transversely the singular spaces $\JJ_{n-1}$.

\begin{prop}\label{splitting} \cite{FHKmem}
For $r<k<n$ there are $\Gamma$-, respectively $\Gamma^\alpha$-equivariant decompositions 
$$\begin{array}{rl}
C_r&=\oplus_{\alpha\in I_r} \big(C_r^\alpha\otimes\Z[\Gamma/\Gamma^\alpha]\big)\\
C_r^\alpha&=\oplus_{\psi\in I_r^\alpha} \big(C_r^{\alpha,\psi}\otimes\Z[\Gamma^\alpha/\Gamma^\psi]\big)\,.
\end{array}$$
Hence, there are homological decompositions
$$\begin{array}{rl}
H_*(\Gamma;C_r)&=\oplus_{\alpha\in I_r} H_*(\Gamma^\alpha;C_r^\alpha)\\
H_*(\Gamma^\alpha;C_r^\alpha)&=\oplus_{\psi\in I_r^\alpha} H_*(\Gamma^\psi;C_r^{\alpha,\psi})\,.
\end{array}$$\qed
\end{prop}

\subsection{A homological framework}
We develop further tools from homological algebra for working with these and associated sequences of modules. A standard background text for this material is Wiebel's book \cite{CW}.

\begin{definition} Let $\mathcal M_*$ be the category of bounded $\Z$-graded $\Z\Gamma$ complexes. Thus an object in $\mathcal M_*$ is a finite sequence of $\Z\Gamma$-modules and maps
$$0\longrightarrow M_s\bra\delta M_{s-1}\bra\delta M_{s-2}\longrightarrow\cdots\longrightarrow M_{t+1}\bra\delta M_t\longrightarrow 0$$
for some $s\geqslant t$ with $\delta^2=0$. Each module is assigned a $\Z$-valued grading, and $\delta$ is a degree $-1$ homomorphism, i.e., reduces grading by 1. Morphisms in $\mathcal M_*$ are degree preserving commutative maps of such complexes. We shall typically denote objects of $\mathcal M_*$ by underlined letters while non-underlined letters are individual $\Z\Gamma$-modules. If $\uM_*\in\mathcal M_*$, denote by $\uM_*[r]$ the complex with the same modules and $\delta$-maps as $\uM_*$, but with degrees increased by $r$, i.e., if $M_s$ occurs in $\uM_*$ in degree $s$, it occurs in $\uM_*[r]$ in degree $s+r$. Unless otherwise stated, a module denoted $M_s$ will be understood to be in degree $s$; in our sequences such as (\ref{complex}), the final copy of $\Z$ is in degree $-1$.
\end{definition}

If $N$ is any individual $\Z\Gamma$-module, we shall at times wish to consider it as an object in $\mathcal M_*$ namely the complex with just one non-zero entry, namely $N$ in degree 0. In the same way we shall write $N[r]$ for the object in $\mathcal M_*$ with just one non-zero entry, namely $N$ in degree $r$.

\smallskip
For $\uM_*\in\mathcal M_*$, denote by $H_*(\uM)$ the homology of the complex $\uM_*$, i.e.,
$$H_r(\uM)=\ker(M_r\bra\delta M_{r-1})/\im(M_{r+1}\bra\delta M_r)\,.$$

\begin{definition} Let $\uM_*\in\mathcal M_*$. Define $H_*(\Gamma;\uM_*)$ as the total homology of the chain complex $\uP_*\otimes_{\Z\Gamma}\uM_*$ where $\uP_*$ is any projective $\Z\Gamma$ resolution of $\Z$. Without loss, we may consider $\uP_*$ to be a free resolution. Recall that if $\uM_*$ and $\uN_*$ are objects in $\mathcal M_*$ with boundary maps $\delta_M$ and $\delta_N$, the total complex of the product $\uM_*\otimes_{\Z\Gamma}\uN_*$ has as module in degree $s$ the sum $\oplus_{p+q=s}M_p\otimes N_q$ and boundary map $\delta_M\otimes 1 + (-1)^p\otimes\delta_N$.
\end{definition} 

Note that $H_s(\Gamma;\uM_*)=H_{s+r}(\Gamma;\uM_*[r])$. 

\smallskip
We also note the standard property that an exact sequence of objects $0\to \underline{A}_*\to \underline{B}_*\to \underline{C}_*\to 0$ in $\mathcal M_*$, i.e., maps of complexes which are exact in each degree, gives rise to a long exact sequence on taking homology $H_*(\Gamma;-)$. (For simplicity we shall denote by $0$ the zero complex in $\mathcal M_*$ consisting of the zero module in every degree.) We also note that, as usual, there are two spectral sequences computing the total homology, one beginning with the double complex $\uP_*\otimes_{\Z\Gamma}\uM_*$ and taking first the homology with respect to the boundary maps in $\uM_*$, the second beginning with $\uP_*\otimes_{\Z\Gamma}\uM_*$ but taking first the homology with respect to the boundary map in the $\Z\Gamma$ resolution $\uP_*$. An immediate consequence of the first of these spectral sequences is the following observation.

\begin{lemma}
If $\uM_*\in\mathcal M_*$ is exact, then $H_*(\Gamma;\uM_*)=0$.\qed
\end{lemma}

\begin{lemma}\label{homequiv}
Suppose 
$$0\longrightarrow M_s\bra\delta M_{s-1}\bra\delta M_{s-2}\longrightarrow\cdots\longrightarrow M_{t+1}\bra\delta M_t\longrightarrow 0$$
is exact, and for some $s\geqslant r> t$, write $\uX_*$ and $\uY_*$ for the complexes
$$\begin{array}{rl}
\uX_*:&0\longrightarrow M_{r-1}\longrightarrow\cdots\longrightarrow M_{t+1}\bra\delta M_{t}\longrightarrow0\,.\\
\uY_*:&0\longrightarrow M_s\bra\delta M_{s-1}\longrightarrow\cdots\longrightarrow M_{r}\longrightarrow0\\
\end{array}$$
Then $H_{i-1}(\Gamma;\uX_*)=H_{i}(\Gamma;\uY_*)=H_{i-r}(\Gamma; K)$ where $K$ is the kernel of the map $\delta\co  M_{r-1}\to M_{r-2}$ (i.e., the image of $M_r\to M_{r-1}$) but considered to be in degree 0.
\end{lemma}

\smallskip\noindent{\bf Proof.} The inclusion and projection maps make $0\to\uX_*\to \underline{M}_*\to \uY_*\to 0$  exact and the left hand equality follows from the induced long exact sequence in group homology and the previous lemma. The right hand equality comes by computing, taking the initial differential that in the graded coefficient module.\qed

\subsection{Exact sequences for pattern cohomology}
We turn now to the specific element of $\mathcal M_*$ we wish to study, namely the exact sequence (\ref{complex}) which for convenience we shall denote $\uC_*$. We define some auxiliary subcomplexes as follows
$$\begin{array}{rl}
\uA_*:&0\to C_{n-1}\to\cdots\to C_0\to 0\,;\\
\uT_*:&0\to C_{n}\to\cdots \to C_0\to 0\,;\\
\uD^r_*:&0\to C_{n}\to\cdots \to C_{r}\to 0\,,\qquad n\geqslant r\geqslant0.\\
\end{array}$$

\begin{lemma}\label{HofP}
There is a $\Gamma$-equivariant equivalence 
$$H_*(\uA)\cong \lim_i H_*(\JJ_{n-1}(i))\,.$$
\end{lemma}

\smallskip\noindent{\bf Proof.} By the lemma \ref{CWcomp} the space $\JJ_{n-1}(i)$ is a CW complex whose $r^{\rm th}$ cellular chain group in the limit as $i\to\infty$ is $C_r$. The complex $\uA_*$ is defined as the cellular chain complex of this space. \qed

\smallskip
As in \cite{FHKcmp, FHKmem} we write $C^0_r$ for ker$(C_{r}\to C_{r-1})$, so there is an exact sequence
\begin{equation}\label{C0seq}
0\to C^0_r\to C_r \to C_{r-1}\to\cdots C_0\to\Z\to 0\,.
\end{equation}

\begin{lemma}\label{identifying}
$$\begin{array}{rcl}
H_*(\Gamma;C_{r-1}^0)&=&H_{*+r}(\Gamma;\uD^{r}_*)\\
H_*(\Gamma;C_n)&=&H_{*+n}(\Gamma;\uD^n_*)\,.
\end{array}$$
\end{lemma}
\smallskip\noindent{\bf Proof.} The first equality follows from (\ref{C0seq}) and lemma \ref{homequiv}. The second from identifying $\uD_*^n$ with $C_n[n]$.\qed

\smallskip
The calculations of \cite{FHKcmp, FHKmem} progressed by inductively working with long exact sequences in group homology given by the short exact sequences of modules
\begin{equation}\label{coeffSESs}\begin{array}{ll}
&0\to C_0^0\to C_0\to\,\Z\,\to0\,,\\
&0\to C_1^0\to C_1\to C_0^0\to0\,,\\
&\cdots \\
&0\to C_n \to C_{n-1} \to C_{n-2}^0\to0\,.
\end{array}
\end{equation}
where the maps $C_q\to C_{q-1}^0$ are induced by the maps $C_q\to C_{q-1}$ and the exactness of $\uC_*$. 

The last of these exact sequences, and one we shall concentrate on later, runs
\begin{equation}\label{FHKLES}
\cdots \to H_{*+1}(\Gamma;C_{n-2}^0)\to H_*(\Gamma;C_n)
 \to H_*(\Gamma;C_{n-1})\to H_*(\Gamma;C_{n-2}^0)\to\cdots\,.
 \end{equation}

\begin{remark}In \cite{FHKmem} these long exact sequences were collected together into a single spectral sequence. From our perspective in this paper, this is the spectral sequence induced by the filtration of $D^n_*=C_n[n]$ given by
\begin{equation}\uT_*=\uD^0_*\to \uD^1_*\to\cdots\to \uD^n_*=C_n[n]\,.\end{equation}
To see the equivalence it is enough to note that the exact sequence of coefficient modules $0\to C^0_r\to C_r\to C_{r-1}^0\to 0$ gives rise to the same long exact sequence in group cohomology as the exact sequence in $\mathcal M_*$
$$0\to C_{r}[r]\to \uD^{r}_*\to \uD^{r+1}_*\to 0$$
though care needs to be taken to check that the degrees and the maps between groups correspond as claimed; we omit the details as the observation is not central to the work which follows.

In particular, however, we note that the long exact sequence of \cite{FHKcmp, FHKmem}, namely (\ref{FHKLES}) above, is induced by the short exact sequence 
\begin{equation}\label{FHKles}0\to C_{n-1}[n-1]\to \uD^{n-1}_*\to C_n[n]\to0\,.\end{equation}
\end{remark}

\begin{remark}\label{mucoeffs}
In the $\mathcal M_*$ framework, the connecting maps $H_s(\Gamma;C^0_{r})\to H_{s-1}(\Gamma;C^0_{r+1})$ in the long exact sequences arising from (\ref{coeffSESs}) correspond to the maps $H_{s+r+1}(\Gamma;\uD_*^{r+1})\to H_{s+r+1}(\Gamma;\uD_*^{r+2})$. Thus  the iterated sequence of connecting maps
$$H_s(\Gamma;\Z)\to H_{s-1}(\Gamma;C_0^0)\to\cdots\to H_{s-n}(\Gamma;C_n)$$
which occurs in our later calculations can be identified with the map in $H_s(\Gamma;-)$ induced by the projection $\uT=\uD_*^0\to \uD^{n}_*=C_n[n]$. 
\end{remark}

The algebraic framework we have set up allows for other exact sequences in homology. In particular, we have the following analogue of Kalugin's sequence \cite{Kal}, though our construction does not need the rationality constructions of \cite{Kal} (in fact, it can be set up without even requiring the earlier assumption that the cut \& project scheme is \can). Consider the short exact sequence in $\mathcal M_*$
\begin{equation}\label{Pses}
0\to \uA_*\buildrel j\over\longrightarrow \uT_*\buildrel m\over\longrightarrow  C_n[n]\to0\,.
\end{equation}
This yields a long exact sequence
\begin{equation}\label{FHKPavel}
\cdots\to H_{*+n}(\Gamma;\uA_*)\buildrel j_*\over\longrightarrow  H_{*+n}(\Gamma;\uT_*)\buildrel m_*\over\longrightarrow  H_{*}(\Gamma;C_n)\to H_{*+n-1}(\Gamma;\uA_*)\to\cdots\,.
\end{equation}
In the next section, under an additional assumption, we will provide a geometric realisation of this sequence, identifying it more explicitly with that of \cite{Kal}. It will relate the Cech cohomology of $\MP$, namely $H^*(\Gamma;C_n)$, with the homology of the $N$-torus $\bT$ given by $H_*(\Gamma;\uT_*)$ and the group homology determined by the complex $\uA_*$, which will be identified with the homology of a certain subspace $\bA$ of $\bT$.

\begin{remark}
The long exact sequence (\ref{FHKPavel}) in fact follows directly from the total homology of the double complex $\uP_*\otimes_{\Z\Gamma}\uA_*$. Computing the total homology by first taking homology with respect to the differential for $\uA_*$ produces an $E^2$-page of the spectral sequence given by
$$E^2_{p,q}=H_p(\Gamma;H_q(\uA))=\left\{
\begin{array}{ll}
H_p(\Gamma;C_n)&\mbox{if $q=n-1$}\\
H_p(\Gamma;\Z)&\mbox{if $q=0$}\\
0&\mbox{otherwise.}
\end{array}\right.$$
The line for $q=0$ is of course the same as $H_p(\Gamma;\uT_*)$ by Lemma \ref{homequiv}. There can only be one more differential, namely $d_n\co  H_*(\Gamma;\uT_*)\to H_{*-n}(\Gamma;C_n)$ and the homology of this computes $H_*(\Gamma;\uA_*)$, giving as it does the long exact sequence (\ref{FHKPavel}).\end{remark}

The following result directly links the two sequences (\ref{FHKLES}) and (\ref{FHKPavel}) and hence the two approaches of \cite{FHKcmp,FHKmem} and \cite{Kal}, a comparison result which will be useful in our computations of $H^*(\MP)$ in the final section. 

\begin{prop}\label{ladder}
There is a commutative diagram 
$$\begin{array}{ccccccccc}
\cdots\to\!\!\!& H_{*+n}(\Gamma;\uA_*)&\!\!\!\buildrel j_*\over\longrightarrow\!\!\!& H_{*+n}(\Gamma;\uT_*)&\!\!\!\buildrel m_*\over\longrightarrow\!\!\!& H_{*}(\Gamma;C_n)&\!\!\!\to\!\!\!& \!\! H_{*+n-1}(\Gamma;\uA_*)\!\!\! &\!\!\!\!\to\cdots\\
&\big\downarrow&&\big\downarrow&&\phantom{cong}\big\downarrow\cong&&\big\downarrow&\\
\cdots \to\!\!& \!\! H_{*+1}(\Gamma; C_{n-1})\!\! &\!\!\!\to\!\!\! &\!\!\! H_{*+1}(\Gamma;C^0_{n-2})\!\!\! 
& \!\!\!\to\!\!\! &H_*(\Gamma;C_n)&\!\!\!\to \!\!\!&H_*(\Gamma;C_{n-1})&\!\!\!\to\cdots\ \
\end{array}$$
in which the rows are exact.
\end{prop}

\smallskip\noindent{\bf Proof.}
The obvious inclusion and projection maps yield the following commutative diagram in which the rows are exact.
$$\begin{array}{ccccccc}
0\to& \uA_*&\buildrel j\over\longrightarrow& \uT_*&\buildrel m\over\longrightarrow& C_n[n]&\to0\\
&\downarrow&&\downarrow&&\vert\!\vert
&\\
0\to& C_{n-1}[n-1]&\to& \uD^{n-1}_*&\to& \uD^n_*
&\to 0\\
\end{array}
$$
On identifying the groups and degrees, this induces the commutative diagram of long exact sequences as in the statement of the Proposition.\qed

\section{Geometric realisation}\label{GeoReal}

In this section we introduce the rationality conditions which allow us to realise various of the elements of $\mathcal M_*$ of the last section and their group homologies in terms of finite cell complexes. This will aid computation in the more difficult examples at the end of the paper. We relate the conditions to the combinatorial condition that the number $L_0$ is finite, equivalently to the condition that the cohomology groups $H^*(\MP)$ are finitely generated.

Assume we have an \can\ cut \& project scheme, and so there is a set $\JJ_{n-1}$ of singular $(n-1)$-dimensional affine subspaces of $\Eperp$, and we have chosen a finite set $\W=\{W_\alpha\}_{\alpha\in I_{n-1}}$ of affine subspaces generating $\JJ_{n-1}$ as $\JJ_{n-1}=\{W_\alpha+\Gammaperp\}_{\alpha\in I_{n-1}}$. Intersections of the elements of $\JJ_{n-1}$ form the lower dimensional singular spaces, or are empty.  Each singular space $U\in\JJ_r$ has associated to it the subgroup $\Gamma^U$ of $\Gamma$ which stabilises $U$ under the natural (projected) translation action of $\Gamma$. 

\begin{definition} A {\em rational subspace\/} of $E$ is a subspace spanned by vectors from $\Q\Gamma$. A {\em rational affine subspace\/} of $E$ is a translate of a rational subspace.
\end{definition}

\begin{definition}\label{ratcond}
A {\em rational projection method pattern\/} is any point pattern arising from an \can\ cut \& project scheme satisfying the following {\em rationality conditions}. 
\begin{enumerate}
\item The number $\nu =\frac{N}{n}= 1+\frac{d}{n}$ is an integer.
\item There is a finite set $\D$ of rational affine subspaces of $E$ in one to one correspondence under $\pi^\perp$ with the set $\W$, i.e., each $W\in\W$ is of the form $W=\pi^\perp(D)$ for some unique $D\in\D$.
\item The members of $\D$ are $\nu(n-1)$-dimensional, and any intersection of finitely many members of $\D$ or their translates is either empty or a rational affine subspace $R$ of dimension $\nu\dim\pi^\perp(R)$.
\end{enumerate}
\end{definition} 


\noindent Extending the notation of Section \ref{Sect2}, for any affine subspace $R$ in $E$, we denote by $\Gamma^R$ the stabiliser subgoup of $\Gamma$ under its translation action on $E$. The following observations are immediate from the geometric set-up.

\begin{lemma}
Suppose we have a rational projection pattern with data as in the definition above. Suppose the singular space $U$ in $\Eperp$ corresponds to some rational affine subspace $R$ in $E$ with $U=\pi^\perp(R)$. Then the stabiliser subgroups of both $U$ and $R$ coincide and the rank of this subgroup equals the dimension of $R$ as an affine subspace.\qed
\end{lemma}

\begin{example}
Consider the {\em Ammann-Beenker}, or {\em Octagonal\/} scheme -- for details see, for example, \cite{Beenker}. In this scheme we have $E=\R^4$ with $\Gamma=\Z^4\subset\R^4$ the integer lattice. Let $v_i$, $i=1,\ldots, 4$ be the four unit vectors 
$$(1,0,0,0)\qquad (0,1,0,0)\qquad (0,0,1,0)\qquad (0,0,0,1)$$
which both generate $\Gamma$ and form a basis for $E$. Consider the linear map $\R^4\to\R^4$ given with respect to this basis by the matrix
$$\left(\begin{array}{cccc}
0&1&0&0\\0&0&1&0\\0&0&0&1\\-1&0&0&0\\
\end{array}\right)\,.$$
This is a rotation of order 8 and has two 2-dimensional eigenplanes, one where the action is rotation by $\pi/4$, the other by $3\pi/4$; take the former for $E^\|$ and the latter for $\Eperp$. Let $W_i$, $i=1,\ldots, 4$, be the 1-dimensional subspace of $\Eperp$ spanned by $\pi^\perp(v_i)$. A set $\W$ generating the singular subspaces is given by $\{W_i\}_{i=1,\ldots, 4}$. The $W_i$ form four rotationally symmetric lines in $\Eperp$ with $W_{i+1}$ the rotation of $W_i$ through $\pi\over 8$. The stabiliser of each $W_i$ is of rank 2: specifically the stabilisers are
$$\begin{array}{rl}
\Gamma^{W_1}=\langle v_1, v_2-v_4\rangle\,,&
\Gamma^{W_2}=\langle v_2, v_1+v_3\rangle\,,\\
\Gamma^{W_3}=\langle v_3, v_2+v_4\rangle\,,&
\Gamma^{W_4}=\langle v_4, v_1-v_3\rangle\,.
\end{array}$$
There is a rational affine plane arrangement $\D=\{D_i\}_{i=1,\ldots, 4}$ covering this family $\W$ where each $D_i$ is the 2-dimensional subspace in $\R^4$ defined by taking as basis the generators of $\Gamma^{W_i}$, as listed above.
\end{example}
\smallskip

\begin{remark}
  Even for \can\ schemes with $\nu$ an integer, it is not always immediately
  clear when there exists a finite set $\D$ of affine planes
  satisfying the rationality conditions. We shall see in Corollary
  \ref{fingen} that if $H^*(\MP)$ is not finitely generated
  (equivalently, if $L_0$ is infinite) then there cannot be a
  lift. However, conversely, suppose $H^*(\MP)$ is finitely generated,
  then Theorem \ref{L0stuff} tells us that the rank of the stabiliser
  $\Gamma^W$ of each dimension $n-1$ singular plane $W\subset \Eperp$
  is $\nu(n-1)$ and so any lift $D$ of $W$ must be an affine space
  parallel to the subspace spanned by the elements of
  $\Gamma^W\subset\Gamma$; the issue is {\em which\/} parallel plane
  to choose, in particular, how to make the relative choices of lifts
  over all the $W\in\W$. 
  
  Along the lines of the discussion at the end
  of the Appendix of \cite{Kal}, in the case where we can choose
  singular planes $\W$ all meeting in a common intersection point, a
  solution is easily given by choosing any point in $E$ over this
  intersection point as a intersection point of the $D\in\D$. This is
  the situation, for example, in the canonical case, where
  $\Gamma=\Z^N$, the integer lattice in $E$, and the acceptance window
  is the $\pi^\perp$-projection of the unit cube, but this is
  certainly not the only situation that allows lifts $\D$. 
  
  Slightly more generally, instead of a common intersection point we can request
  that each $W\in\W$ contains some rational point with respect to a
  basis of $\Gamma^\perp$ and a suitably chosen origin. We say then that
  $W$ has also rational position, in addition to the rational orientation.
  Since intersections of affine spaces in rational position and orientation
  also have rational position and orientation, all singular spaces then
  have rational positions. In fact, such a singular space in rational
  position and orientation contains a dense subset of rational points, 
  a rational affine subspace, whose rational dimension is equal to
  the rank of the stabilizer in $\Gamma$. Such rational affine spaces 
  have a preferred lift with the required properties. We choose as 
  origin of $E$ a point above the origin of $E^\perp$, and as lattice 
  basis of $\Gamma$ the unique lift $\{b_i\}$ of the chosen basis 
  $\{b^\perp_i\}$ of $\Gamma^\perp$. Every rational point  
  $\sum_iq_ib^\perp_i$ is then lifted to $\sum_iq_ib_i$, and rational 
  affine subspaces of $E^\perp$ are thus lifted to rational affine 
  subspaces of $E$ of the same rational dimension. As the full lift
  of a singular subspace we thus take the closure of the lift of its
  rational subset. With this scheme, the lift of the intersection of 
  two affine subspaces is always equal to the intersection of the two 
  lifts, as required. 
  
  The situation with singular spaces in rational 
  position actually includes the case with a common intersection point 
  of all $W\in\W$, but is still by no means the most general one. 
  The generalised Penrose patterns \cite{PenGen} are examples of
  rational projection patterns where the elements of $\W$ have 
  positions which can move continuously when the parameter $\gamma$ 
  is varied, and which do not have a common intersection point.
\end{remark}


Given a rational projection scheme, denote by $\RR_{n-1}$ the set $\D+\Gamma$ of all $\nu(n-1)=(N-\nu)$-dimensional affine subspaces in the $\Gamma$ orbit of $\D$. In Section \ref{HomAlg} it was useful to view $\JJ_{n-1}$, the set of all singular $(n-1)$-spaces in $\Eperp$, as the increasing union of locally finite collections of $(n-1)$-spaces, $\JJ_{n-1}=\cup_i\JJ_{n-1}(i)$. In the same way, denote by $\RR_{n-1}(i)$ the $\nu (n-1)$-dimensional affine subspaces which correspond to the elements of $\JJ_{n-1}(i)$. Again for convenience, we also denote by $\RR_{n-1}$ and $\RR_{n-1}(i)$ the subspaces of $E$ consisting of the union of the subspaces in these sets.


\begin{lemma}
The projection map $\pi^\perp$ induces homology isomorphisms 
$$H_*(\RR_{n-1}(i))\cong H_*(\JJ_{n-1}(i))\,.$$
\end{lemma}

\smallskip\noindent{\bf Proof.}
The homologies $H_*(\RR_{n-1}(i))$ and $H_*(\JJ_{n-1}(i))$ may each be computed, in principle, by Mayer-Vietoris spectral sequences corresponding to the construction of $\RR_{n-1}(i)$ and $\JJ_{n-1}(i)$ as unions of $\nu(n-1)$-  and $(n-1)$-dimensional planes respectively. The map $\pi^\perp$ induces a one-to-one correspondence between the planes and intersection planes in $\RR_{n-1}(i))$ and $\JJ_{n-1}(i)$, and as in both cases each such plane is contractible, $\pi^\perp_*$ induces an isomorphism on the first page of the spectral sequence, and hence an isomorphism of the final homologies.\qed

\begin{cor}\label{georeal}
The projection map $\pi^\perp$ induces $\Gamma$-equivariant isomorphisms 
$$H_*(\RR_{n-1})\cong H_*(\uA)\qquad\mbox{and}\qquad
H_*(E)\cong H_*(\uT)\,.$$
\end{cor}

\smallskip\noindent{\bf Proof.} For the first, as $\RR_{n-1}$ is a CW complex and can be considered as the direct limit of the $\RR_{n-1}(i)$ we have $H_*(\RR_{n-1})\cong\lim_i H_*(\RR_{n-1}(i))$ since homology commutes with direct limits. The previous lemma gives an equivalence
$$\lim_i H_*(\RR_{n-1}(i))\cong\lim_i H_*(\JJ_{n-1}(i))$$
and the right hand object is equivalent to $H_*(\uA)$ by Lemma \ref{HofP}. The second equivalence is immediate since $\Eperp$ is a ($\Gamma$-equivariant) homotopy retract of $E$.\qed

\begin{definition}
Write $\bA$ for the quotient space $\RR_{n-1}/\Gamma$ and $\bT$ for $E/\Gamma$. Write $\alpha$ for the induced inclusion $\alpha\co \bA\to\bT$. Clearly $\bT$ is just the $N$-torus.
\end{definition}

\begin{theorem}\label{identification}
There is a commutative diagram whose vertical maps are isomorphisms
$$\begin{array}{ccc}
H_*(\bA)&\buildrel\alpha_*\over\longrightarrow &H_*(\bT)\\
\vert\!\vert&&\vert\!\vert\\
H_*(\Gamma;\uA_*)&\buildrel j_*\over\longrightarrow &H_*(\Gamma;\uT_*)
\end{array}$$
and $j_*$ is induced by the inclusion $\uA_*\to\uT_*$ as in the exact sequence (\ref{Pses}).
\end{theorem}

\smallskip\noindent{\bf Proof.}
The quotient maps $\RR_{n-1}\to\RR_{n-1}/\Gamma=\bA$ and $E\to E/\Gamma=\bT$ induce fibrations and maps
$$\begin{array}{cccccc}
\RR_{n-1}&\to&\bA&\to&B\Gamma\\
\big\downarrow&&\phantom{\alpha}\big\downarrow\alpha&&\big\vert\!\big\vert\\
E&\to&\bT&\to&B\Gamma
\end{array}$$
where $B\Gamma$ is the classifying space of the group $\Gamma$. These lead to computations of $H_*(\bA)$ and $H_*(\bT)$ via Serre spectral sequences, which compute these homologies as the total homologies of the double complexes $\uP_*\otimes_{\Z\Gamma} C_*(\RR_{n-1})$ and $\uP_*\otimes_{\Z\Gamma} C_*(E)$ where $\uP_*$ as in Section 3 is any free $\Z\Gamma$ resolution of $\Z$ while $C_*(\RR_{n-1})$ and $C_*(E)$ are $\Gamma$-chain complexes computing the homologies of $\RR_{n-1}$ and $E$ respectively.

By Corollary \ref{georeal}, after the first differential of the spectral sequences, the resulting double complexes are identical to those computing respectively $H_*(\Gamma;\uA_*)$ and $H_*(\Gamma;\uT_*)$. Moreover, the map of double complexes induced by $\alpha$ is from this point on identical to that induced by $j\co  \uA_*\to\uT_*$.\qed

\smallskip
The long exact sequence (\ref{FHKPavel}) may now be interpreted as follows, recovering the exact sequence of \cite{Kal}. For simplicity, we denote the homomorphism $H_*(\bT)\to H_*(\Gamma;C_n[n])$ given by the composite of $m_*$ with the identification $H_*(\Gamma;\uT_*)\cong H_*(\bT)$ of Theorem \ref{identification} also by $m_*$.

\begin{cor}\label{Ples} There is an exact sequence
$$\begin{array}{ccccccccc}
\cdots \to &H_r(\bA)&\buildrel\alpha_*\over\longrightarrow &H_r(\bT)&
\buildrel m_*\over\longrightarrow &H_r(\Gamma;C_n[n])&\longrightarrow &H_{r-1}(\bA)&\to\cdots\\
&&&&&\vert\!\vert&&&\\
&&&&&H^{N-r}(\MP)&&&\\
\end{array}$$\qed
\end{cor}

\begin{remark}\label{identMU}
Strictly speaking, to fully identify this sequence with that of \cite{Kal} we need to show that the composite $H_r(\bT)\to H_{r-n}(\Gamma;C_n)\cong H^{N-r}(\MP)$ can be identified with the map $\mu^*\co  H^{N-r}(\bT)\to H^{N-r}(\MP)$ composed with the Poincare duality isomorphism $H_r(\bT)\cong H^{N-r}(\bT)$. This can be done by identifying the action of $\mu^*$ with the map in group cohomology $H^*(\Gamma;-)$ induced by the coefficient map $\uT_*\to C_n[n]$ as in Remark \ref{mucoeffs}. We briefly return to this issue in the Appendix, as the complete identification requires the construction realising $H^*(\MP)$ as the group cohomology $H^*(\Gamma;C_n)$, but for now we omit the details as this point is not necessary for the work which follows.
\end{remark}

As $\bA$ is a cell complex with top cells of dimension $(N-\nu)$, we have $H_r(\bA)=0$ for $r>N-\nu$. Corollary \ref{Ples} immediately gives

\begin{cor}\label{PThigh} For a rational projection pattern, there are isomorphisms $H_r(\Gamma;\Z)\cong H_r(\Gamma;C_n)$ for $r>N-\nu+1$. Equivalently, there are isomorphisms $H^s(\bT)\cong H^s(\MP)$ for $s<\nu-1$. \qed
\end{cor}

\begin{cor}\label{fingen}
For any commutative ring $S$, the cohomology $H^*(\MP;S)$ of a rational projection pattern $P$ is finitely generated over $S$.
\end{cor}

\smallskip\noindent{\bf Proof.} Recall that if $X$ is a space with the homotopy type of a finite CW complex, then $H_*(X)$ is finitely generated over $\Z$. The spaces $\bA$ and $\bT$ both have the homotopy type of finite CW complexes, and hence so too has the mapping cone $C(\alpha)$ of $\alpha\co \bA\to\bT$. The exact sequence of Corollary \ref{Ples} says that $H^{N-*}(\MP)\cong H_*(C(\alpha))$ and hence the groups are finitely generated. The result for general $S$ follows by a standard universal coefficient theorem argument.\qed

\medskip
The advantage of Theorem \ref{identification} is that it allows information useful for computing with the long exact sequences (\ref{FHKles}), (\ref{FHKPavel}) to be obtained from the reasonably tangible map of topological spaces $\bA\to\bT$. The subspace $\bA$ of the torus $\bT$ is itself given as the union of $(N-\nu)$-tori, each such torus being $T_i=D_i/\Gamma^{D_i}$ as $D_i$, $i\in I_{n-1}$, runs over the elements of $\D$. A consequence of the rationality conditions means that any intersection of finitely many of these tori is either empty or a common subtorus of the form $R/\Gamma^R$. This structure, together with details of the data describing the rational affine subspaces $\D$, makes $H_*(\bA)$ and the homomorphism $\alpha_*$ accessible, at least in principle: for any given projection scheme of course, the finite complex $\bA$ can of course have considerable complexity. 

The following observations specify the main phenomena that specific computation must address. Rewriting the exact sequence of Corollary \ref{Ples}, we obtain
\begin{equation}\label{Basicalpha}
0\to \coker(\alpha_*)\to H^*(\MP)\to \ker(\alpha_*)\to0\,.
\end{equation}
Thus for computations in rational cohomology, it suffices to compute the ranks of the homomorphisms $\alpha_*\co H_r(\bA;\Q)\to H_r(\bT;\Q)$. However, for {\em integral\/} computations, there are potential extension problems to solve if there is torsion in $\ker(\alpha_*)$, which will certainly be the case if there is torsion in $H_*(\bA)$, since $H_*(\bT)$ is torsion free.

As noted in the proof of Corollary \ref{fingen}, there is an isomorphism 
$$H^{N-r}(\MP)\cong H_r(C(\alpha))\,,$$ where $C(\alpha)$ is the mapping cone of the map $\alpha$ (equivalently, $H_*(C(\alpha))$ is the relative homology $H_*(\bT,\bA)$). Given the finite generation result, Corollary \ref{fingen}, we know by the universal coefficient theorem (UCT)  going between homology and cohomology the groups $H_*(C(\alpha))$ (and hence $H^*(\MP)$) if we can compute the {\em cohomology\/} $H^*(C(\alpha))$. Explicitly, and as regards torsion components, the torsion subgroup of $H^{N-r}(\MP)$, which is the torsion subgroup of $H_r(C(\alpha))$, is isomorphic to the torsion subgroup of $H^{r+1}(C(\alpha))$. This latter cohomology group sits in an extension analogous to (\ref{Basicalpha}) (i.e., the long exact sequence in cohomology of the pair $(\bT,\bA)$)
\begin{equation}\label{COHses}
\begin{array}{rl}0\to\coker\!\!\!\!\!\!&\left(H^r(\bT)\buildrel\alpha^*\over\longrightarrow H^r(\bA)\right)
\to H^{r+1}(C(\alpha))\\&\qquad\qquad\qquad\to
\ker\left(H^{r+1}(\bT)\buildrel\alpha^*\over\longrightarrow H^{r+1}(\bA)\right)\to0\,.
\end{array}\end{equation}
Note that the right hand group, the $\ker$-term, is here necessarily torsion free, since $H^*(\bT)$ is. Thus this short exact sequence splits and the only torsion component in $H^*(C(\alpha))$ must arise as the torsion component of the $\coker$-term of (\ref{COHses}). Explicitly, let us define
$$\begin{array}{rl}
s_r\ =&\mbox{free abelian rank of }\ker\left(H^{r}(\bT)\buildrel\alpha^*\over\longrightarrow H^{r}(\bA)\right)\\
f_r\ =&\mbox{free abelian rank of }\coker\left(H^{r}(\bT)\buildrel\alpha^*\over\longrightarrow H^{r}(\bA)\right)\\
\mathcal T_r\ =&\mbox{torsion subgroup of }\coker\left(H^r(\bT)\buildrel\alpha^*\over\longrightarrow H^r(\bA)\right)\,.
\end{array}$$
Then $H^{r+1}(C(\alpha))=\Z^{f_{r}+s_{r+1}}\oplus \mathcal T_r$, and by the UCT
\begin{cor}\label{idtorsion}
The cohomology group $H^{N-r}(\MP)=H_r(C(\alpha))$ is given by
$$H^{N-r}(\MP)=\Z^{f_{r-1}+s_{r}}\oplus \mathcal T_r\,.$$\qed
\end{cor}

\section{Patterns of codimension one and two}\label{Sect5}

The exact sequence (\ref{COHses}) and Corollary \ref{idtorsion} show that in principle the cohomology groups $H^*(\MP)$ for a rational projection pattern are completely determined by knowledge of the homomorphisms $\alpha^*\co  H^*(\bT)\to H^*(\bA)$. The homology or cohomology of $\bA$  is potentially accessible via a Meyer-Vietoris spectral sequence computation arising from the decomposition of $\bA$ into its component $(N-\nu)$-tori; this is the approach of the calculations (with $\Q$ coefficients) of \cite{Kal}, and we utilise aspects of this approach for some of our work in the final section. 

In this section however, and for our initial work on codimension 3 patterns in Section \ref{Sect6}, we use instead the exact sequence (\ref{FHKLES}) as our fundamental tool and compute inductively up the values of $n$, the codimension. The two approaches are essentially equivalent for rational patterns, but the inductive approach has some merits in terms of spreading out the computations into manageable steps, and in particular is also applicable to patterns not satisfying the rationality conditions. In general, for whichever approach, the complexity and subtlety of the computations increases significantly as $n$ increases. 


\subsection{Codimension 1}
We consider \can\ projection patterns of codimension 1, and note that the faces of a one-dimensional acceptance domain are points and so $\JJ_0$ consists of a finite number of distinct $\Gamma$-orbits of points; as before, $L_0$ denotes the number of these orbits. 

\begin{theorem} 
For a dimension $d$, codimension 1 \can\ projection pattern, $H^{d-k}(\MP)=H_k(\Gamma;C_1)$
is a free abelian group of rank
\begin{equation*}\left\{\begin{array}{ll}
0 & \mbox{for}\quad k>d,\\
{{d+1}\choose{k+1}} & \mbox{for}\quad d\geqslant k>0,\\
L_0+d & \mbox{for}\quad k=0.
\end{array}\right.\end{equation*}

\end{theorem}

\smallskip\noindent{\em Proof.}
We compute  $H_*(\Gamma;C_1)$ using the short exact sequence of 
$\Z\Gamma$-modules
$$0\longrightarrow C_1\longrightarrow
C_0\longrightarrow\Z\longrightarrow 0\,,$$ 
which is the complex (\ref{complex}) for $n=1$. In this sequence, $\Z$
carries the trivial $\Gamma$ action, while the action of $\Gamma$ on
$C_0$ is free. In group homology we get the long exact sequence
$$\cdots\to H_{k+1}(\Gamma;\Z)\to H_k(\Gamma;C_1)\to
H_k(\Gamma;C_0)\to\cdots\,.$$ 
Now $H_k(\Gamma;\Z)\cong\Lambda_k\Gamma$ 
is just the homology of a $(d+1)$-torus, so
$H_k(\Gamma;\Z)$ is free abelian of rank 
${d+1}\choose{k}$. 
Meanwhile, the freeness of $\Gamma$ on $C_0$ means that the homology
groups $H_k(\Gamma;C_0)$ 
are zero for $k>0$ and $H_0(\Gamma;C_0)=\Z^{L_0}$.

Our long exact sequence now tells us  that
$H_k(\Gamma;C_1)\cong\Lambda_{k+1}\Gamma\cong
\Z^{{{d+1}\choose{k+1}}}$ in dimensions $k>0$ and for  dimension 0 there is an exact
sequence
$$0\to\Lambda_1\Gamma \to H_0(\Gamma;C_1)\to
\Z^{L_0}\stackrel{\epsilon}{\to}\Z\to0\,.$$
Hence $H_0(\Gamma;C_1)\cong\Lambda_1\Gamma \oplus \ker \epsilon$,
and so is free abelian of rank
$L_0+d$. \null\hfill$\square$

\begin{remark}
While not needed for the work below, we note in passing that the same result holds for the case where we would allow the acceptance domain to have infinitely many connected components and where $L_0$ may be infinite.
The explicit details needed can be found in \cite{FHKmem} Chapter III where a different approach to the codimension 1 case is taken, and it is shown that $\MP$ can be modelled by 
a punctured torus. 
\end{remark}

\subsection{Codimension 2}\label{cd2}
We turn to the case $n=2$. The theorem below is stated for any \can\ projection pattern with finitely generated cohomology, so in particular holds for any rational projection pattern. Our analysis proceeds via the pair of exact sequences of (\ref{coeffSESs}), 
\begin{equation}\label{2ses}
0\to C_2\to C_1\to
  C^0_0\to0\qquad0\to C^0_0\to C_0\to\Z\to0\,.
\end{equation}
Denote by $\beta_k$ the homomorphism in $H_k(\Gamma;-)$ induced by the module homomorphism $C_1\to C^0_0$; the relevant part of the sequence (\ref{FHKLES}) now runs
\begin{equation}\label{kercoker}
0\to\mbox{coker }\beta_{k+1}\to H_k(\Gamma;C_2)\to\mbox{ker
}\beta_k\to0\,.
\end{equation}
Set $R_k$ to be the rank of $\erz{\Lambda_{k+1}\Gamma^{\alpha}:\alpha\in I_1}$, the subgroup of $\Lambda_{k+1}\Gamma$ generated by all the images of the inclusions $\Lambda_{k+1}\Gamma^\alpha\to\Lambda_{k+1}\Gamma$.

\begin{theorem}\label{thm51}
Let $P$ be an \can\  projection pattern with codimension 2 and suppose $H^*(\MP)$ is finitely generated. Thus, in particular, the dimension $d$ is $2(\nu-1)$ and the numbers $L_1$, $L_0$ and $L_0^\alpha$ are finite. Each group $H^{d-k}(\MP)=H_k(\Gamma;C_2)$ is thus a sum of a free abelian group and a finite abelian torsion group. 
\begin{enumerate}
\item Sequence (\ref{kercoker}) splits and $H_k(\Gamma;C_2)\cong\coker\beta_{k+1}\oplus \ker \beta_k$.\\
\item The rank of the free abelian part of $H_k(\Gamma;C_2) $ is given by the formulae
$$\begin{array} {cl}
 \left( 2\nu\atop 2+k\right)+L_1 \left( \nu\atop 1+k\right) - R_k - R_{k+1},& \mbox{ for }0<k\leqslant d,\mbox{ and}\\ 
 \sum_{j=0}^2 (-1)^j \left( 2\nu\atop 2-j\right) + L_1\sum_{j=0}^1 (-1)^j
\left(\nu\atop 1-j\right) + e  - R_1,& \mbox{ for }k=0,
\end{array}$$ 
where $e$ is  the Euler characteristic and is given by
$$ e=\sum_{p} (-1)^p \rk_\Q H_{p}(\Gamma;C_2) =-L_0+\sum_{\alpha\in I_1}L^\alpha_0 .$$\\
\item The torsion part of $H_k(\Gamma;C_2)$ is given by
the torsion part of the cokernel of $\beta_{k+1}$, which can be identified here as the map
$$\bigoplus_{\alpha\in I_1}\Lambda_{k+2}
\Gamma^\alpha\to\Lambda_{k+2}\Gamma$$  
induced by the inclusions $\Gamma^\alpha\to\Gamma$. 
In particular, $H^{d-k}(\MP)=H_k(\Gamma;C_2)$ is torsion free for
$k\geqslant d/2$.
\end{enumerate}
\end{theorem}

\smallskip\noindent{\bf Proof.}
The right hand sequence of (\ref{2ses}) in group homology behaves identically to the calculations in the previous subsection for codimension 1. We obtain 
$$H_k(\Gamma;C_0^0) \cong
\left\{\begin{array}{ll}
\Z^{\left(d+2 \atop
        k+1\right)}=\Lambda_{k+1}\Gamma & \mbox{for}
\quad
    k>0,\\
\Z^{d+L_0+1} &
    \mbox{for}\quad
k=0
\end{array}\right.
$$
where the $k=0$ case
arises from the short exact sequence
$$0\to\Lambda_1\Gamma\to
H_0(\Gamma;C_0^0)\to \ker\epsilon\to
0\qquad\mbox{ with
}\epsilon\co \Z^{L_0}\to\Z.$$

Using the splitting of Proposition \ref{splitting}, which here identifies
$H_k(\Gamma;C_1)$ with $\bigoplus_\alpha H_k(\Gamma^\alpha;C^\alpha_1)$, a
similar calculation based on the exact sequences
\begin{equation}\label{sesal}
0\to C_1^\alpha \to
  C_0^\alpha \to\Z\to 0
\end{equation}
gives
\begin{equation}\label{co2-H1}
H_k(\Gamma;C_1)
\cong \left\{\begin{array}{ll}
\bigoplus_{\alpha\in
      I_1}\Lambda_{k+1}\Gamma^\alpha & 
\mbox{for}\quad
    k>0,\\
\bigoplus_{\alpha\in I_1}(\Lambda_{1}\Gamma^\alpha\oplus
    
\ker\epsilon^\alpha) 
& \mbox{for}\quad
    k=0
\end{array}\right.
\end{equation}
where $\epsilon^\alpha$ denotes the augmentation $H_0(\Gamma^\alpha;C_0^\alpha)\cong
\Z^{L_0^\alpha}
\to \Z$.
Recall that the rank of each $\Gamma^\alpha$ is $\nu$. The internal direct sum in the case $k=0$ represents the splitting of the short exact sequences
$$0\to\Lambda_1\Gamma^\alpha\to
H_0(\Gamma^\alpha;C^\alpha_1)\to
\ker\epsilon^\alpha\to
0\,.$$

For $k>0$, the homomorphism in $\beta_k\co  H_k(\Gamma;C_1)\to H_k(\Gamma;C^0_0)$ identifies with the  homomorphism
\begin{equation}\label{eq-beta}
\bigoplus_{\alpha\in
  I_1}\Lambda_{k+1}\Gamma^\alpha\to\Lambda_{k+1}\Gamma
  \end{equation}
induced by the inclusions $\Gamma^\alpha\to\Gamma$. Similarly, 
$\beta_0$
identifies with the homomorphism of
extensions
\begin{equation}\label{beta0}\begin{array}{rcccl}
0\to&
\bigoplus_{\alpha\in
      I_1}\Lambda_1\Gamma^\alpha&
\to
    H_0(\Gamma;C_1)\to&
\bigoplus_{\alpha\in
      I_1}\ker\epsilon^\alpha&
\to 0
    
\\
&\qquad\bigg\downarrow\beta_0'&\quad\bigg\downarrow\beta_0&\qquad\bigg\downarrow\beta_0''&
\\
0\to&
\Lambda_1\Gamma&
\to H_0(\Gamma;C^0_0)\to&
\ker\epsilon&
\to 0 
\end{array}\end{equation}
induced by the inclusions
$\Gamma^\alpha\to\Gamma$ and $C_0^\alpha\subset C_0$.

We can now prove the claims of the theorem by organising the data from these calculations and exact sequences. The reader may find it helpful to consult the diagram (\ref{diagramn=2}) which displays this information for the case $\nu=2$; the analogue for higher values of $\nu$ is very similar, though obviously longer in the vertical direction.

For part (1), note that the left hand sequence in (\ref{2ses}) gives the long exact sequence
in homology
$$\cdots\to H_{k+1} (\Gamma;C_1)\buildrel \beta_{k+1}\over\longrightarrow H_{k+1} (\Gamma;C^0_0)\to H_k(\Gamma;C_2)\to H_k(\Gamma;C_1)\buildrel \beta_k\over\longrightarrow\cdots$$
which at $H_k(\Gamma;C^2)$  may be written as the short exact sequence
$$0\to\mbox{coker }\beta_{k+1}\to H_k(\Gamma;C^2)\to\mbox{ker
}\beta_k\to0\,.$$ 
This splits since $\mbox{ker }\beta_k\subset H_k(\Gamma;C_1)$ is finitely generated free abelian.

For part (2) it is sufficient to work with rational coefficients and count ranks. Note that $R_k$ is the rank of the image of $\beta_k$ and that $\beta_0$ is surjective.

For part (3), the torsion part of $H_k(\Gamma;C_2)$ must arise from $\mbox{coker }\beta_{k+1}$ since $\mbox{ker }\beta_k$ is free. However,  as the rank of $\Lambda_{k+2}\Gamma^\alpha$ is ${\nu\choose{k+2}}$, for $k\geqslant d/2=\nu-1$ this is trivial and so in this range the map $\beta_{k+1}$ is zero and there is no torsion in its cokernel.\qed

\smallskip
This final result, putting bounds on where torsion may appear, will be
seen to be a special case of a result for 
arbitrary codimension in Subsection~\ref{Gen}. Examples suggest that
these bounds are best possible. 
\medskip

A direct computation of $\ker\beta_0$ from the data encoded in
$\ker\beta_0'$ and $\ker\beta_0''$ and the diagram (\ref{beta0}) need not be
immediate, a point which will become a serious issue when we deal
with the codimension 3 patterns later. The diagram (\ref{beta0}) gives a
an exact sequence
$$0\to\ker\beta_0'
\to\ker\beta_0
\to\ker\beta_0''
\stackrel{\Delta_0}\to\coker\beta_0'
\to0$$
($\coker\beta_0=0$ as $\beta_0$ is surjective). 
In general there is no reason why the connecting map
$\Delta_0:\ker\beta_0''\to\coker\beta_0'$ should be trivial. In fact,
the T\"ubingen Triangle Tiling \cite{BKSZ,KSB} tiling is an example in
which $\coker\beta'_0=\Z_5$ and hence $\Delta_0$ is
non-trivial. However, we do not need to compute $\Delta_0$ explicitly
as the $\coker\beta'_0$ term will only be comprised of torsion terms,
which do not contribute either to the torsion or the free rank
of the cohomology of the tiling. 

\vfill
\begin{diag}\label{diagramn=2}
The entire computation for the case $\nu=2$, i.e., codimension $=$ dimension $= 2$
can be summarized in the following diagram in which all rows and columns
are exact.

$$
\begin{array}{rcccl}
&0=\bigoplus_{\alpha\in
I_1}\Lambda_4\Gamma^\alpha=\!\!\!\!\!\!\!\!\!\!&H_3(\Gamma;C_1)&&
\\
&&\quad\bigg\downarrow\beta_3&&
\\
&\Z=\Lambda_4\Gamma=\!\!\!\!\!\!\!\!\!\!\!\!\!\!\!\!\!\!\!\!\!\!\!
&H_3(\Gamma;C^0_0)&&
\\
&&\bigg\downarrow&&
\\
&&H_2(\Gamma;C_2)&\!\!\!\!\!\!\!\!\!\!\!\!\!\!\!\!\!\!\!\!\!\!\!\!\!\!\!\!\!\!
=H^0(\MP)&
\\
&&\bigg\downarrow&&
\\
&0=\bigoplus_{\alpha\in
I_1}\Lambda_3\Gamma^\alpha=\!\!\!\!\!\!\!\!\!\!&H_2(\Gamma;C_1)&&
\\
&&\quad\bigg\downarrow\beta_2&&
\\
&\Z^4=\Lambda_3\Gamma=\!\!\!\!\!\!\!\!\!\!\!\!\!\!\!\!\!\!\!\!\!\!\!
&H_2(\Gamma;C^0_0)&&
\\
&&\bigg\downarrow&&
\\
&&H_1(\Gamma;C_2)&\!\!\!\!\!\!\!\!\!\!\!\!\!\!\!\!\!\!\!\!\!\!\!\!\!\!\!\!\!\!
=H^1(\MP)&
\\
&&\bigg\downarrow&&
\\
&\Z^{L_1}=
\bigoplus_{\alpha\in I_1}\Lambda_2\Gamma^\alpha=&
H_1(\Gamma;C_1)&
\\
&
&\quad\bigg\downarrow\beta_1&&
\\
&\qquad\qquad \Z^{{2\nu\choose 2}}=
\Lambda_2\Gamma=&H_1(\Gamma;C^0_0)&
\\
&&\bigg\downarrow&&
\\
&&H_0(\Gamma;C_2)&\!\!\!\!\!\!\!\!\!\!\!\!\!\!\!\!\!\!\!\!\!\!\!\!\!\!\!\!\!\!
=H^2(\MP)&
\\
&&\bigg\downarrow&&
\\
0\to\!\!\!\!\!\!\!\!\!\!\!\!\!\!\!\!\!\!\!\!\!\!&
\bigoplus_{\alpha\in I_1}\Lambda_1\Gamma^\alpha\!\!\!\!\!\!\!\!\!\!\!\!\!\!\!\!\!\!\!\!\!\!&
\to\qquad H_0(\Gamma;C_1)\qquad\to&
\bigoplus_{\alpha\in I_1}\ker\epsilon^\alpha&
\to 0 
\\
&\qquad\qquad\qquad\bigg\downarrow\beta_0'&\quad\bigg\downarrow\beta_0&\qquad\bigg\downarrow\beta_0''&
\\
0\to\!\!\!\!\!\!\!\!\!\!\!\!\!\!\!\!\!\!\!\!\!\!&
\Lambda_1\Gamma\!\!\!\!\!\!\!\!\!\!\!\!\!\!\!\!\!\!\!\!\!\!\!\!\!\!\!\!\!\!&
\to\qquad H_0(\Gamma;C^0_0)\qquad\to&
\ker\epsilon&
\to 0 
\\
&&\bigg\downarrow&&
\\
&&0&&
\end{array}
$$
\end{diag}

\subsection{Example computations}
\label{sec:excomp}

All examples discussed below have to some 
extent been calculated by computer. For this purpose, we have used
the computer algebra system GAP \cite{GAP}, the GAP package Cryst
\cite{EGN97,Cryst}, as well as further software written in the GAP 
language. It should be emphasized that these computations are not 
numerical, but use integers and rationals of unlimited size or 
precision. Neglecting the possibility of programming errors, 
they must be regarded as exact. 

One piece of information that needs to be computed is the set of
all intersections of singular affine subspaces, along with
their incidence relations. This is done with code based on the
Wyckoff position routines from the Cryst package. The set of
singular affine subspaces is invariant under the action of a
space group. Cryst contains routines to compute intersections
of such affine subspaces and provides an action of space group
elements on affine subspaces, which allows to compute space group
orbits. These routines, or variants thereof, are used to determine 
the space group orbits of representatives of the singular affine 
subspaces, and to decompose them into translation orbits. The 
intersections of the affine subspaces from two translation orbits 
is the union of finitely many translation orbits of other affine 
subspaces. These intersections can be determined essentially by 
solving a linear system of equations modulo lattice vectors,
or modulo integers when working in a suitable basis. With these 
routines, it is possible to generate from a space group and a 
finite set $\mathcal{W}$ of representative singular affine spaces 
the set of all singular spaces, their intersections, and their 
incidences.

A further task is the computation of ranks, intersections, and 
quotients of free $\Z$-modules, and of homomorphisms between such 
modules, including their kernels and cokernels. These are standard 
algorithmic problems, which can be reduced to the computation of 
Smith and Hermite normal forms of integer matrices, including the 
necessary unimodular transformations \cite{Coh95}. GAP already 
provides such routines, which are extensively used.

The codimension 2 examples discussed here all have
dihedral symmetry of order $2n$, with $n$ even. The lattice 
$\Gamma^\perp_{n}$ is given by the $\Z$-span of the vectors in the 
star $e_{i}=(\cos(\frac{2\pi i}{n}),\sin(\frac{2\pi i}{n}))$,
$i=0,\ldots,n-1$. The singular lines have special orientations 
with respect to this lattice. They are parallel to mirror lines 
of the dihedral group, which means that they are either {\em along}
the basis vectors $e_{i}$, or {\em between} two neighboring basis 
vectors, i.e., along $e_{i}+e_{i+1}$. In all our examples below, 
with the single exception of the  generalized Penrose tilings \cite{PenGen}, 
one line from each translation orbit passes through the origin. We denote 
the sets of representative singular lines by $\mathcal{W}_{n}^{a}$
and $\mathcal{W}_{n}^{b}$, for lines along and between the basis 
vectors $e_{i}$. The defining data of several well-known tilings
can now be given as a pair of a (projected) lattice, and a set
of translation orbit representatives of singular lines. Specifically, 
the Penrose tiling \cite{deBruijn} is defined by the pair 
$(\Gamma^\perp_{10},\mathcal{W}_{10}^{a})$, the T\"ubingen Triangle 
Tiling (TTT) \cite{BKSZ,KSB} by the pair 
$(\Gamma^\perp_{10},\mathcal{W}_{10}^{b})$, 
the undecorated octagonal Ammann-Beenker tiling \cite{Beenker} by the pair 
$(\Gamma^\perp_{8},\mathcal{W}_{8}^{a})$, and the undecorated Socolar
tiling \cite{Soc89} by the pair $(\Gamma^\perp_{12},\mathcal{W}_{12}^{a})$. 
For the decorated versions of the Ammann-Beenker \cite{Soc89,AGS,G93} 
and Socolar tilings \cite{Soc89}, the set of singular lines 
$\mathcal{W}_{n}^{a}$ has to be replaced by 
$\mathcal{W}_{n}^{a}\cup\mathcal{W}_{n}^{b}$, $n=8$ and $12$, respectively. 
These well-known examples are complemented by the coloured Ammann-Beenker
tiling with data $(\Gamma^\perp_{8},\mathcal{W}_{8}^{b})$, which can be realised
geometrically by colouring the even and odd vertices of the classical
Ammann-Beenker tiling with two different colours, and the heptagonal
tiling from \cite{GK}, which is given by the pair 
$(\Gamma^\perp_{14},\mathcal{W}_{14}^{a})$. Except for the coloured Ammann-Beenker 
tiling, which had not been considered in the literature previously, 
the rational ranks of the cohomology of these tilings had been computed 
in \cite{GK}. Table~\ref{tab2d} shows their cohomology with integer 
coefficients, including the torsion parts where present.

The generalized Penrose tilings \cite{PenGen} are somewhat different 
from the tilings discussed above. They are built upon the decagonal
lattice $\Gamma^\perp_{10}$ too, but have only fivefold rotational symmetry.
The singular lines do not pass through the origin in general,
and their positions depend on a continuous parameter $\gamma$.
For instance, the representative lines of the two translation 
orbits of lines parallel to $e_{0}$ pass through the points 
$-\gamma e_{1}$ and $\gamma(e_{1}+e_{2})$. It turns out that
these shifts of line positions always lead to the same line 
intersections and incidences. Even multiple intersection points
remain stable, and are only moved around if $\gamma$ is varied.
Consequently, all generalized Penrose tilings have the same
cohomology, except for $\gamma\in\Z[\tau]$, which corresponds
to the real Penrose tilings \cite{deBruijn}. This had already been 
observed by Kalugin \cite{Kal}, and is in contradiction with 
the results given in \cite{GK}, which were obtained due to a 
wrong parametrisation of the singular line positions. Corrected 
results are given in Table~\ref{tab2d}.

\begin{table}
\caption{Cohomology of codimension 2 tilings with dihedral symmetry.}
\label{tab2d}
\begin{center}
\small
\newcommand{\op}{\oplus}
\renewcommand{\arraystretch}{1.3}
\begin{tabular}{|l|c|c|c|}
\hline
Tiling                       & $H^2$                &$H^1$      & $H^0$  \\
\hline
Ammann-Beenker (undecorated) & $\Z^{ 9}$            & $\Z^{ 5}$ & $\Z^1$ \\
Ammann-Beenker (coloured)     & $\Z^{14}\op\Z_2$     & $\Z^{ 5}$ & $\Z^1$ \\
Ammann-Beenker (decorated)   & $\Z^{23}$            & $\Z^{ 8}$ & $\Z^1$ \\
Penrose                      & $\Z^{ 8}$            & $\Z^{ 5}$ & $\Z^1$ \\
generalized Penrose          & $\Z^{34}$            & $\Z^{10}$ & $\Z^1$ \\
T\"ubingen Triangle          & $\Z^{24}\op\Z_5^{2}$ & $\Z^{ 5}$ & $\Z^1$ \\
Socolar (undecorated)        & $\Z^{28}$            & $\Z^{ 7}$ & $\Z^1$ \\
Socolar (decorated)          & $\Z^{59}$            & $\Z^{12}$ & $\Z^1$ \\
\hline
\end{tabular}
\end{center}
\end{table}

Among the tilings discussed above, only the TTT and the coloured
Ammann-Beenker tiling have torsion in their cohomology. 
The set of singular lines of the TTT is constructed from the 
lines $\mathcal{W}_{10}^{b}$. The translation stabilizers 
$\Gamma^{\alpha}$ of all these lines are contained in a common 
sublattice $\Gamma'^\perp_{10}$ generated by the star of vectors 
$e_{i}+e_{i+1}$; it has index 5 in $\Gamma^\perp_{10}$.
It is therefore not too surprising that $\coker\beta_{1}$
(Theorem~\ref{thm51}) develops a torsion component $\Z_{5}^{2}$, 
which shows up in the cohomology group $H^2$ of the TTT, in 
agreement with the results obtained using the method of Anderson-Putnam
\cite{AP} which computes the cohomology of TTT via its substitution structure.
In much the same way, and for analogous reasons, a torsion component 
$\Z_{2}$ in $H^2$ is obtained also for the coloured Ammann-Beenker tiling,
and also the four-dimensional, codimension 2 tilings with 
data $(\Gamma^\perp_{14},\mathcal{W}_{14}^{b})$ have torsion 
components $\Z_{7}^{4}$ in $H^4$, and $\Z_{7}^{3}$ in $H^{3}$,
in agreement with the bounds given in Theorem~\ref{thm51}.

There is an interesting relation between the TTT and the Penrose
tiling. Since the lattice $\Gamma'^\perp_{10}$ is rotated by $\pi/10$ 
with respect to $\Gamma^\perp_{10}$, the TTT can also be constructed
from the pair $(\Gamma'^\perp_{10},\mathcal{W}_{10}^{a})$. However,
the singular set $\Gamma'^\perp_{10}+\mathcal{W}_{10}^{a}$ is even
invariant under all translations from $\Gamma^\perp_{10}$, so that it
is equal to $\Gamma^\perp_{10}+\mathcal{W}_{10}^{a}$, which defines
the Penrose tiling. In other words, the TTT and the Penrose tiling 
have the same set of singular lines, only the lattice $\Gamma^\perp$ 
acting on it is different. The TTT is obtained by breaking the 
translation symmetry of the Penrose tiling to a sublattice of 
index 5. This explains why the Penrose tiling is locally derivable 
from the TTT, but local derivability does not hold in the opposite
direction \cite{BSJ91}. A broken symmetry can be restored in a local way, 
but the full lattice symmetry cannot be broken to a sublattice 
in any local way, because there are no local means to distinguish 
the five cosets of the sublattice. Any tiling whose set of singular 
lines accidentally has a larger translation symmetry are likely 
candidates for having torsion in their cohomology.

For the coloured Ammann-Beenker tiling, the situation is completely
analogous. Geometrically, the coloured and the uncoloured version
are the same, and thus have the same singular lines $\mathcal{W}_{8}^{a}$.
For the coloured variant, we have to restrict the lattice to the
colour preserving translations, which form a sublattice $\Gamma'^\perp_8$
of index 2 in $\Gamma^\perp_8$. With respect to $\Gamma'^\perp_8$, the
lines in $\mathcal{W}_{8}^{a}$ are between the generating vectors, 
so that the pair $(\Gamma'^\perp_8,\mathcal{W}_{8}^{a})$ is equivalent 
to the pair $(\Gamma^\perp_8,\mathcal{W}_{8}^{b})$. Again, the uncoloured
Ammann-Beenker tiling can be recovered from the coloured one by restoring
the translations broken by the colouring.

\section{Patterns of codimension three}\label{Sect6}

The case of projection patterns of codimension 3 is a good deal more
complex than the codimension 2 theory, though the principles of 
computation remain the same. We shall initially consider \can\ projection patterns with finitely generated cohomology (equivalently, that $L_0$ is finite). Later we shall specialise to the rational projection patterns.

The dimension 3 space $\Eperp$ now has families of singular lines and
singular planes. Following \cite{FHKmem} we shall index by 
$\theta$ the lines and by $\alpha$ the planes. The rank of the main
group $\Gamma$ is $N$ and, by Theorem \ref{L0stuff}, the rank of the stabiliser $\Gamma^\alpha$ of a
singular plane is $N-\nu=2\nu$, while the stabiliser $\Gamma^\theta$ of a singular 
line has rank $\nu$.

The complex (\ref{complex}) in this case can be broken into two exact sequences
\begin{equation}\label{n=3complex}
0\to C_3\to C_2\to C_1^0\to0\quad\mbox{and}\quad 0\to C_1^0\to C_1\to C_0\to\Z\to0\,.
\end{equation}
The left hand sequence gives the long exact sequence
\begin{equation*}
\cdots\to H_{s}(\Gamma;C_3) \to H_{s}(\Gamma;C_2)
\buildrel\phi_{s}\over\longrightarrow H_{s}(\Gamma;C_1^0)\to
H_{s-1}(\Gamma;C_3)\to\cdots
\end{equation*}
computing $H_*(\Gamma;C_3)=H^{d-*}(\MP)$ so long as we know the groups and homomorphisms
$$\phi_s\co  H_s(\Gamma;C_2)\to H_s(\Gamma;C_1^0)$$
and can solve the resulting extension problems. The groups $H_s(\Gamma;C_1^0)$ are computed from the right hand sequence of (\ref{n=3complex}) following exactly the same procedure we used to compute the codimension 2 examples from the analogous complex. We obtain

\begin{lemma}\label{HC10} There are equalities and short exact sequences
\begin{equation*}
\begin{array}{rll}
H_s(\Gamma;C_1^0)&=0 &\mbox{for}\quad s\geqslant N-1,\\
H_s(\Gamma;C_1^0)&=\Lambda_{s+2}\Gamma&\mbox{for}\quad N-1>s\geqslant \nu,\\
0\to\Lambda_{s+2}\Gamma\to H_s(\Gamma;C_1^0)&\to\ker\gamma_s\to0\, ,&\mbox{for}\quad \nu>s\geqslant 1,\\
0\to\coker\gamma_1\to H_0(\Gamma;C_1^0)&\to\ker\gamma_0\to0\, ,&
\end{array}
\end{equation*}
where, for $s>0$,
$$\gamma_s\co  H_s(\Gamma;C_1)=\bigoplus_{\theta\in
  I_1}\Lambda_{s+1}\Gamma^\theta\to
H_s(\Gamma;C^0_0)=\Lambda_{s+1}\Gamma\, ,$$ 
and
$$\begin{array}{rcccl}
0\to&
\bigoplus_{\theta\in I_1}\Lambda_1\Gamma^\theta&
\to H_0(\Gamma;C_1)\to&
\bigoplus_{\theta\in I_1}\ker\epsilon^\theta&
\to 0 
\\
&\qquad\bigg\downarrow\gamma_0'&\quad\bigg\downarrow\gamma_0&
\qquad\bigg\downarrow\gamma_0''&
\\
0\to&
\Lambda_1\Gamma&
\to H_0(\Gamma;C^0_0)\to&
\ker\epsilon&
\to 0 
\end{array}$$
are both induced by the inclusions $\Gamma^\theta\to\Gamma$ and
$C^\theta_0\subset C_0$. 
Note that all the terms in (\ref{HC10}) are free of torsion except
possibly the coker$\,\gamma_1$ summand.  
\end{lemma}

By Proposition \ref{splitting} the groups $H_*(\Gamma;C_2)$ split as $\bigoplus_{\alpha\in I_2} H_*(\Gamma^\alpha;C_2^\alpha)$ and for each singular plane $C_2^\alpha$ we have a sequence
$$0\to C_2^\alpha\to C_1^\alpha\to C_0^\alpha\to\,\Z\,\to0\,.$$
As before, we obtain 
\begin{lemma}\label{HC2}
\begin{equation*}
\begin{array}{rll}
H_s(\Gamma^\alpha;C_2^\alpha)&=0 &\mbox{for}\quad
s\geqslant 2\nu-1,\\
H_s(\Gamma^\alpha;C_2^\alpha)&=\Lambda_{s+2}\Gamma^\alpha&\mbox{for}\quad
2\nu-1>s\geqslant \nu,\\ 
0\to\Lambda_{s+2}\Gamma^\alpha\to
H_s(\Gamma^\alpha;C_2^\alpha)&\to\ker\beta^\alpha_s\to0\, ,&\nu>s\geqslant 1\\
0\to\coker\beta^\alpha_1\to
H_0(\Gamma^\alpha;C_2^\alpha)&\to\ker\beta^\alpha_0\to0\, ,&
\end{array}
\end{equation*}
where, for $s>0$,
$$\beta^\alpha_s\co  H_s(\Gamma;C_1^\alpha)=\bigoplus_{\theta\in
  I_1^\alpha}\Lambda_{s+1}\Gamma^\theta\to
H_s(\Gamma^\alpha;C_0^{\alpha0})=\Lambda_{s+1}\Gamma^\alpha\, ,
$$
and
\begin{equation}\begin{array}{rcccl}
0\to&
\bigoplus_{\theta\in I_1^\alpha}\Lambda_1\Gamma^\theta&
\to H_0(\Gamma;C^\alpha_1)\to&
\bigoplus_{\theta\in I_1^\alpha}\ker\epsilon^\theta&
\to 0 
\\
&\qquad\bigg\downarrow\beta^\alpha_0\null'&\quad\bigg\downarrow\beta^\alpha_0&
\qquad\bigg\downarrow\beta^\alpha_0\null'\null'&
\\
0\to&
\Lambda_1\Gamma^\alpha&
\to H_0(\Gamma;C_0^{\alpha0})\to&
\ker\epsilon^\alpha&
\to 0 .
\end{array}\nonumber\end{equation}
The maps $\beta_s^\alpha$ are again induced by the obvious inclusions, and the
only potential torsion term in (\ref{HC2}) arises from the
$\coker\beta^\alpha_1$ expression;  all other terms are free
abelian.
\end{lemma}

The expression for $\phi_s\co  H_s(\Gamma;C_2)\to H_s(\Gamma;C_1^0)$
under the identifications (\ref{HC10},$\,$\ref{HC2}) can be
obtained as in \cite{FHKmem}:
$\phi_s$ is a 
sum of morphisms $\phi_s^\alpha$, which in turn are determined by the
diagrams 
\begin{equation}\label{eq-phi}
\begin{array}{rcccccl}
{\to} & H_s(\Gamma;C^\alpha_2\otimes\Z[\Gamma/\Gamma^\alpha]) & \to &
\bigoplus_{\theta\in I_1^\alpha} \Lambda_{s+1}\Gamma^{\theta} &
\bra{\beta^\alpha_s} & \Lambda_{s+1}\Gamma^{\alpha} & {\to}\\
&\bigg\downarrow \phi_s^\alpha & & \bigg\downarrow
j^\alpha_{s} & & \bigg\downarrow \imath^\alpha_s& \\ {\to} &
H_s(\Gamma;C_1^0) & \to & \bigoplus_{\theta\in I_1}
\Lambda_{s+1}\Gamma^{\theta} & \bra{\gamma_s} &
\Lambda_{s+1}\Gamma & {\to} 
\end{array}
\end{equation}
for $s>0$, and
\begin{equation*}
\begin{array}{rcccccl}
{\to}&\!\!\!H_0(\Gamma;C^\alpha_2\otimes\Z[\Gamma/\Gamma^\alpha])\!\!\! &
\to &
\!\!\!\bigoplus_{\theta\in I_1^\alpha}(\Lambda_{1}\Gamma^{\theta}
\oplus\ker\epsilon^\theta)\!\!\!  &
\bra{\beta^\alpha_0} & \!\!\!\Lambda_{1}\Gamma^{\alpha}
\oplus\ker\epsilon^\alpha\!\!\!& 
{\to0}\\
&\bigg\downarrow \phi_0^\alpha & & \bigg\downarrow
j^\alpha_{0} & & \bigg\downarrow \imath^\alpha_0& \\ {\to} &
H_0(\Gamma;C_1^0) & \to & \!\!\!\bigoplus_{\theta\in I_1}
(\Lambda_{1} \Gamma^{\theta}\oplus\ker\epsilon^\theta)\!\!\! &
\bra{\gamma_0} &
\!\!\!\Lambda_{1}\Gamma\oplus\ker\epsilon\!\!\! & {\to0} 
\end{array}
\end{equation*}
where $j_s^\alpha$ and $\imath^\alpha_s$ are induced by the obvious
inclusions.

\begin{diag}\label{diagramn=3}
The entire computation for the case $\nu=2$ can be summarized in the following diagram in which all rows and columns are exact. The general case is similar but with a longer diagram. The finite generation condition on cohomology means that the only  rational projection patterns in dimension 3 are those with codimension 1 or 3.
\vfill

\begin{equation}\label{dia1}
\begin{array}{rcccl}
&0=\bigoplus_{\alpha\in
I_2}\Lambda_6\Gamma^\alpha=\!\!\!\!\!\!\!\!\!\!&H_4(\Gamma;C_2)&&
\\
&&\bigg\downarrow&&
\\
&\Z=\Lambda_6\Gamma=\!\!\!\!\!\!\!\!\!\!\!\!\!\!\!\!\!\!\!\!\!\!\!
&H_4(\Gamma;C_1^0)&&
\\
&&\bigg\downarrow&&
\\
&&H_3(\Gamma;C_3)&\!\!\!\!\!\!\!\!\!\!\!\!\!\!\!\!\!\!\!\!\!\!\!\!\!\!\!\!\!\!
=H^0(\MP)&\!\!\!\!\!\!\!\!\!\!\!\!\!\!\!\!\!\!\!\!\!\!\!=\,\Z
\\
&&\bigg\downarrow&&
\\
&0=\bigoplus_{\alpha\in
I_2}\Lambda_5\Gamma^\alpha=\!\!\!\!\!\!\!\!\!\!&H_3(\Gamma;C_2)&&
\\
&&\bigg\downarrow&&
\\
&\Z^6=\Lambda_5\Gamma=\!\!\!\!\!\!\!\!\!\!\!\!\!\!\!\!\!\!\!\!\!\!\!
&H_3(\Gamma;C_1^0)&&
\\
&&\bigg\downarrow&&
\\
&&H_2(\Gamma;C_3)&\!\!\!\!\!\!\!\!\!\!\!\!\!\!\!\!\!\!\!\!\!\!\!\!\!\!\!\!\!\!
=H^1(\MP)&
\\
&&\bigg\downarrow&&
\\
0\to&
\bigoplus_{\alpha\in I_2}\Lambda_4\Gamma^\alpha&
\to H_2(\Gamma;C_2)\to&
\!\!\!\!\!\!\!\!\!\!\!\!\!\!\!\!\!\!\!\!\!\!\!\!\!\!\!\!\!\!0& 
\\
&\bigg\downarrow\phi_2'&\quad\bigg\downarrow\phi_2&&
\\
0\to&
\Lambda_4\Gamma&
\to H_2(\Gamma;C_1^0)\to&
\!\!\!\!\!\!\!\!\!\!\!\!\!\!\!\!\!\!\!\!\!\!\!\!\!\!\!\!\!\!0& 
\\
&&\bigg\downarrow&&
\\
&&H_1(\Gamma;C_3)&\!\!\!\!\!\!\!\!\!\!\!\!\!\!\!\!\!\!\!\!\!\!\!\!\!\!\!\!\!\!
=H^2(\MP)&
\\
&&\bigg\downarrow&&
\\
0\to&
\bigoplus_{\alpha\in I_2}\Lambda_3\Gamma^\alpha&
\to H_1(\Gamma;C_2)\to&
\bigoplus_{\alpha\in I_2}\ker\beta_1^\alpha&
\to 0 
\\
&\bigg\downarrow\phi_1'&\quad\bigg\downarrow\phi_1&\bigg\downarrow\phi_1''&
\\
0\to&
\Lambda_3\Gamma&
\to H_1(\Gamma;C_1^0)\to&
\ker\gamma_1&
\to 0 
\\
&&\bigg\downarrow&&
\\
&&H_0(\Gamma;C_3)&\!\!\!\!\!\!\!\!\!\!\!\!\!\!\!\!\!\!\!\!\!\!\!\!\!\!\!\!\!\!
=H^3(\MP)&
\\
&&\bigg\downarrow&&
\\
0\to&
\bigoplus_{\alpha\in
I_2}\coker\beta_1^\alpha\!\!\!\!\!&
\to H_0(\Gamma;C_2)\to&
\bigoplus_{\alpha\in I_2}\ker\beta_0^\alpha&
\to 0 
\\ 
&\bigg\downarrow\phi_0'&\quad\bigg\downarrow\phi_0&\bigg\downarrow\phi_0''&
\\ 
0\to&
\coker\gamma_1&
\to H_0(\Gamma;C_1^0)\to&
\ker\gamma_0&
\to 0 
\\
&&\bigg\downarrow&&
\\
&&0&&
\end{array}
\end{equation}
\end{diag}

As for the case $n=2$ we shall compute first with rational coefficients and in so doing compute the ranks of the free abelian part of the integral cohomology $H^*(\MP)$, and second consider the torsion part. Although the rational computation amounts to counting dimensions and using the extension
$$0\to \coker\phi_{s+1}\to H_s(\Gamma; C_3)\to\ker\phi_s\to0$$
the computation in terms of accessible numbers does not follow immediately by chasing Diagram \ref{diagramn=3} or its analogue for higher $\nu$: for example, the rank of $\coker\phi_{s+1}$ is not automatically the sum of the ranks of $\coker\phi_{s+1}'$ and $\coker\phi_{s+1}''$, and likewise for the kernels. As before, a simple application of the snake lemma tells us that there are six term exact sequences
\begin{equation}\label{LESdelta}
0\!\to\!\ker\phi_s'
\!\to\!\ker\phi_s
\!\to\!\ker\phi_s''
\!\buildrel\Delta_s\over\longrightarrow\!\coker\phi_s'
\!\to\!\coker\phi_s
\!\to\!\coker\phi_s''
\!\to\!0.\\
\end{equation}
As in the codimension $2$ case, direct knowledge of $\Delta_0$ is unnecessary to solve for either the rational ranks or the torsion. However, the  maps $\Delta_s$, $s>0$, enter into consideration in both cases; of course $\Delta_s$ is trivial for $s\geqslant\nu$ since $\ker\phi_s''=0$ in these degrees. 

The following lemma gives a useful link between the data required for computations using the set-up just described, based on the long exact sequence (\ref{FHKLES}), and the approach to computing $H^*(\MP)$ which uses the sequence (\ref{FHKPavel}). Note that in the case of a rational projection pattern, $j_*$ is essentially the homomorphism $\alpha_*$ of Corollary \ref{Ples}. It will give us a helpful criterion for deciding when $\Delta_s$ vanishes. 

\begin{lemma}\label{cokernels}
For $s>0$, the cokernel of  $j_*\co  H_{s+2}(\Gamma;\uA_*\otimes R)\to H_{s+2}(\Gamma;\uT_*\otimes R)$ is identical to $\coker(\phi'_s)/\im\Delta_s$. In particular, $\Delta_s=0$ if and only if $\coker(j_*)=\coker(\phi'_s)$.
\end{lemma}

\smallskip\noindent{\bf Proof.}   By the exact sequence (\ref{FHKPavel}), $\coker (j_*)=\im (m_*)$. As we can write $m$ as the composite  $\uT_*\to\uD^2_*\to\uD^3_*=C_3[3]$ in $\mathcal M_*$,
we may identify the composite 
$$\Lambda_{s+2}\Gamma\to H_s(\Gamma;C_1^0)\to H_{s+2}(\Gamma;C_3[3])=H_{s-1}(\Gamma;C_3)$$ 
in Diagram \ref{diagramn=3} (or its analogue for higher values of $\nu$) with $m_*$ in $H_*(\Gamma;-)$. Thus $\im(m_*)=\coker (j_*)$ can be identified with the image of the composite
$$\coker(\phi_s')\to\coker(\phi_s)\hookrightarrow H_{s-1}(\Gamma;C_3)$$
which, by the exact sequence (\ref{LESdelta}), is equal to $\coker(\phi'_s)/\im\Delta_s$.\qed




\subsection{Rational computations}

\smallskip
The computation of the rational ranks, i.e., $\dim H_*(\Gamma;C_3\otimes\Q)$, for an \can\ projection pattern is now essentially straightforward, albeit longwinded. Clearly 
$$\mbox{dim }H_s(\Gamma;C_3\otimes\Q)=\mbox{ dim }\coker\phi_{s+1}\,+\,\mbox{ dim }\ker\phi_s$$
and the computations follow from knowledge of the dimensions of the groups $H_*(\Gamma;C_2\otimes\Q)$ and $H_*(\Gamma;C_1^0\otimes\Q)$ obtained via the methods for $n=2$, together with a computation of the ranks of the maps $\phi_s$. The latter are obtained relatively straightforwardly for $s\geqslant\nu$, but for smaller values of $s$ the terms arising via $\phi_s''$ add a further degree of complexity and require knowledge of $\Delta_s$. The following summarises the computation in the general case, and also corrects an error in the determination of the kernel of $\gamma$ in \cite{FHKmem} (middle of page 112). Applied to
the Ammann-Kramer tiling  \cite{AmmannKramer}, these formulae evaluate to agree with the results of \cite{Kal}. All of the terms  used in the statement can be calculated relatively easily on a computer provided $\Delta_s=0$, which we shall see below is the case for any rational projection tiling when using rational coefficients.

\begin{theorem}\label{RatFormulae}{\em [Erratum\footnote{The formulae given in \cite{FHKmem} are correct only if the equation $\langle\im j^\alpha_s\cap\ker\gamma_s:\alpha\in I_2\rangle=\ker\gamma_s$ in the middle of page 112 holds and $\rk\Delta_s=0$. For the Ammann-Kramer tiling and the dual canonical $D_6$ tiling this is not so and the rank of the left hand side is one lower than that of the right.} to Theorem~2.7 of
\cite{FHKmem}] }\label{thm-corrections}
Given an \can\ projection pattern with $L_0$ finite, codimension $3$ and dimension $d=3(\nu-1)$, the following formulae give the ranks of the rational homology groups $H_*(\Gamma;C_3\otimes\Q)$. All ranks are understood to be rational ranks. For $s>0$, 
\begin{eqnarray*} 
\rk H_s(\Gamma;C^n\otimes\Q) &=& \left( 3\nu\atop s+3\right)
+L_2 \left( 2\nu\atop s+2\right) + \sum_{\alpha\in I_2}L_1^\alpha
\left( \nu\atop s+1\right)\\
&&  + L_1\left( \nu\atop s+2\right) - R_s - R_{s+1}, \\ 
\rk H_0(\Gamma;C^n\otimes\Q) & =& \sum_{j=0}^3 (-1)^j
\left( 3\nu\atop 3-j\right)+ L_2\sum_{j=0}^2 (-1)^j \left(
2\nu\atop 2-j\right) \\
&& + \sum_{\alpha\in I_2}L_1^\alpha \sum_{j=0}^1 (-1)^j
\left(\nu\atop 1-j\right) + L_1\sum_{j=0}^2 (-1)^j
\left(\nu\atop 2-j\right)\\
&& + e  - R_1 \, .
\end{eqnarray*} 
Here $R_s = \rk \phi_s + \sum_{\alpha\in I_2}\rk \beta_s^\alpha-\rk \gamma_s+\rk\Delta_s$ which is given by, for $s>1$ 
\begin{eqnarray*}
R_s & = & \rk\erz{\Lambda_{s+2}\Gamma^{\alpha}:\alpha\in I_2}
+\sum_{\alpha\in I_2}\rk\erz{\Lambda_{s+1}\Gamma^{\theta}:\theta\in I_1^\alpha}\\
&& -\rk\erz{\Lambda_{s+1}\Gamma^{\theta}:\theta\in I_1} + \rk\Delta_s
\end{eqnarray*}
and
\begin{eqnarray*}
R_1 &=& \rk\erz{\Lambda_{3}\Gamma^{\alpha}/\im\beta^\alpha:\alpha\in I_2}
+\sum_{\alpha\in I_2}
\rk\erz{\Lambda_{2}\Gamma^{\theta}:\theta\in I_1^\alpha}\\
& &+\rk\erz{(\bigoplus_{\theta\in I_1^\alpha}\Lambda_{2}\Gamma^\theta)\cap
  \ker\gamma_1 : \alpha\in I_2}-\rk\erz{\Lambda_{2}\Gamma^{\theta}:\theta\in I_1} + \rk\Delta_1\, .
\end{eqnarray*}
Finally, the Euler characteristic $ e:=\sum_{s} (-1)^s \rk_\Q H_{s}(\Gamma;C^n)$ is given by 
$$ e =
L_0-\sum_{\alpha\in I_2} L^\alpha_0 +\sum_{\alpha\in
I_2}\sum_{\theta\in I_1^\alpha}L^\theta_0 -\sum_{\theta\in
I_1}L^\theta_0. $$
\end{theorem}

\begin{remark}
For $s\geqslant\nu$ the expression for $R_s$
simplifies to $R_s=  \rk\erz{\Lambda_{s+2}\Gamma^{\alpha}:\alpha\in I_2}$. This follows from the fact that
$\Lambda_{s+1}\Gamma^{\theta}$ vanishes for $\theta\in I_1$ or $I_1^\alpha$ as
the rank of $\Gamma^{\theta}$ is $\nu$. For $\nu=2$, that is if the dimension is 3, the expression for $R_1$ also simplifies slightly as $\im\beta_2^\alpha$ vanishes for similar reasons.
\end{remark}

Now assume the projection pattern considered satisfies the rationality conditions, and so we can use the geometric realisation $\bA\buildrel\alpha\over\longrightarrow\bT$ of the homomorphism $\uA_*\buildrel j\over\longrightarrow\uT_*$ as in Theorem \ref{identification}. In particular, Lemma \ref{cokernels} tells us  that $\Delta_s=0$ if $\im(\alpha_*)=\im(\phi_s')$.

\begin{lemma}\label{FullalphaQ}
For a rational projection pattern, in computations of cohomology with any coefficient ring $R$ a field of characteristic 0, the homomorphisms $\Delta_s=0$ for all $s>0$.
\end{lemma}

\smallskip\noindent{\bf Proof.} Recall that $\bA$ is given as a union of $(N-\nu)$-tori $T_i$ inside $\bT$. The individual inclusions $T_i\to\bT$ combine to give a map factoring  $\coprod T_i\to\bA\buildrel\alpha\over\longrightarrow\bT$ which shows that we always have  the inclusion $\im(\phi_s')\subset\im(\alpha_*)$; we prove the opposite inclusion. It is sufficient to work with the field $\F=\R$. 

Consider a simplicial decomposition of the pair $(\bT, \bA)$, that is, a simplicial decomposition of $\bT$ such that each (open) cell has either empty intersection with $\bA$ or is contained in it. The map on simplicial chain groups given by mapping the simplex $(x_0,\ldots, x_r)$ to $(x_1-x_0)\wedge\cdots\wedge (x_r-x_0)\in\Lambda_r\R\Gamma$ vanishes on boundaries and hence induces an isomorphism between $H_r(\bT;\R)$ and $\Lambda_r\Gamma\otimes\R$. Restricting to $\bA=\cup T_i$ it follows that $\im(\alpha_*)$ is contained in the subgroup of $\Lambda_r\Gamma\otimes\R$ generated by the subgroups $\Lambda_r\Gamma^{D_i}\otimes\R$.\qed

\begin{cor}\label{Ratsimplify}
For rational projection tilings $\rk\Delta_s=0$ for all $s>0$ and the formulae in the statement of Theorem \ref{RatFormulae} correspondingly simplify.\qed
\end{cor}

\subsection{Torsion and the integral computations}
We turn to the determination of the integral cohomology of \can\ projection patterns and, given the results above for calculations with rational coefficients, this entails an examination of the torsion groups which can arise in the computations, and the solution of associated extension problems. As before, we assume throughout this subsection that the number $L_0$ is finite, but do not as yet assume the rationality conditions.

The results of Lemmas \ref{HC10} and \ref{HC2} show that computations for $H_s(\Gamma;C_3)$ are relatively straightforward for $s\geqslant\nu$; in these cases we have an extension
$$\begin{array}{rl}
0\to\coker\big(\phi'_{s+1}\co \bigoplus_{\alpha\in I_2}&\!\!\!\!\!\Lambda_{s+3}\Gamma^\alpha\to\Lambda_{s+3}\Gamma\big)
\to H_s(\Gamma;C_3)\\ \to&
\ker\big((\phi'_{s}\co \bigoplus_{\alpha\in I_2}\Lambda_{s+2}\Gamma^\alpha\to\Lambda_{s+2}\Gamma\big)\to0\\
\end{array}$$
and this extension splits since the kernel term is free abelian. In fact as the rank of $\Gamma^\alpha$ is $2\nu$ we immediately recover the result of Corollary \ref{PThigh} for codimension 3 patterns. As the cokernel term can in principle have torsion, this same observation about the rank of $\Gamma^\alpha$ gives the following analogue of the final line in Theorem \ref{thm51}; again examples suggest this result is best possible.

\begin{prop}\label{notor}
For a codimension 3 \can\ projection pattern, there is no torsion in $H_s(\Gamma;C_3)=H^{d-s}(\MP)$  for $s\geqslant 2(\nu-1)$.\qed
\end{prop}

For the remainder of the paper we specialise to the case $\nu=2$, whose details are depicted in Diagram \ref{diagramn=3}. The general case is similar, but with analogous extension problems arising over a larger range of dimensions. The following summarises the situation and follows immediately from the previous observations.

\begin{theorem}\label{codim3}
For a codimension 3, dimension 3 \can\ projection pattern with finitely generated cohomology (i.e., $L_0$ finite), we have
$$H^{3-s}(\MP)=H_s(\Gamma;C_3)=\left\{\begin{array}{ll}
0& \mbox{for}\quad s\geqslant4,\\
\Z& \mbox{for}\quad s=3,\\
\Z^6\oplus\ker\left\{\phi_2'\co \bigoplus_{\alpha\in I_2}
(\Lambda_4\Gamma^\alpha)\to\Lambda_4\Gamma\right\} 
&\mbox{for}\quad s=2,
\end{array}\right.$$
and so there is no torsion in these degrees. 

In homological degree 1 there is no torsion in $\ker\phi_1$ and we obtain
$$H^2(\MP)=H_1(\Gamma;C_3)=\coker\phi_2'\oplus\ker\phi_1\,.$$
The summand $\coker\phi_2'$ may contain torsion, but this is computable
from the description of $\phi_2'$ as the homomorphism 
$\bigoplus_{\alpha\in I_2}
(\Lambda_3\Gamma^\alpha)\to\Lambda_3\Gamma$ induced by the inclusions
$\Gamma^\alpha\to\Gamma$. 
\end{theorem}

In homological degree 0 (i.e., computing $H^3(\MP)$) the computational problems are considerable. Using this approach alone, we can only deduce the
homology group as the extension
\begin{equation}\label{extn}0\to\coker\phi_1\to H_0(\Gamma;C^3)\to\ker\phi_0\to0\,.\end{equation}
Here torsion can arise in both $\coker\phi_1$ and $\ker\phi_0$;
note that it is not necessarily the case that even if there is torsion in $\ker\phi_0$ then it lifts to torsion elements in $H_0(\Gamma;C^3)$: any specific calculation therefore needs to determine the torsion in  $\coker\phi_1$ and whether any torsion of $\ker\phi_0$ lifts to $H_0(\Gamma;C^3)$. Torsion in $\ker\phi_0$ can arise only from torsion in $\ker\phi_0'$ by Lemma \ref{HC2}.

However, even if $\Delta_1=0$, neither is the determination of torsion in  $\coker\phi_1$ straightforward. Torsion elements in this group may arise from either $\coker\phi_1'$ or $\coker\phi_1''$, but torsion in $\coker\phi_1''$ itself does not immediately imply that it lifts to torsion in $\coker\phi_1$, and there is another extension problem to solve on the way. 

\smallskip
In general such extension problems need further geometric or topological input to solve. To that end, we shall now assume that our pattern satisfies the rationality conditions of Section \ref{GeoReal}. As the cohomology of a rational projection pattern is always finitely generated, Corollary \ref{fingen}, the free abelian part of $H^*(\MP)$ is completely determined by Theorem \ref{RatFormulae} and Corollary \ref{Ratsimplify}. The torsion is described, as in Corollary \ref{idtorsion}, as the torsion subgroup of $\coker\left(H^3(\bT)\buildrel\alpha^*\over\longrightarrow H^3(\bA)\right)$, or alternatively via an extension, as in exact sequence (\ref{Basicalpha}), of the torsion in $H_2(\bA)$ with $\coker(\alpha_*)$. We consider first the computation of the homology and cohomology of the space $\bA$.

We start with the homology of $\bA$. For elementary reasons, $H_0(\bA)=\Z$ and $H_1(\bA)=H_1(\bT)=\Z^6$. To compute the higher homology groups we use the Mayer-Vietoris spectral sequence for the homology of the resolution space $\bA^\Delta$. This considers $\bA$ as the union of $L_2$ 4-tori, $L_1$ 2-tori and $L_0$ 0-tori (points). The $E^1$-page of the spectral sequence has, as its $r^{\rm th}$ column $E^1_{r,*}$, the homology of the disjoint union of those tori arising as $(r+1)$-fold intersections, as shown in Table \ref{tab:MVSS}.

\begin{table*}
\begin{center}
{
\renewcommand{\arraystretch}{1.4}
\begin{tabular}{|c|c|c|}
\hline
$\bigoplus_{\alpha\in I_{2}}\Lambda_{4}\Gamma^{\alpha}$ & & \\
\hline
$\bigoplus_{\alpha\in I_{2}}\Lambda_{3}\Gamma^{\alpha}$ & & \\
\hline
$\bigoplus_{\alpha\in I_{2}}\Lambda_{2}\Gamma^{\alpha}\oplus
 \bigoplus_{\theta\in I_{1}}\Lambda_{2}\Gamma^{\theta}$
  & $\bigoplus_{\alpha\in I_{2}}\bigoplus_{\theta\in I_{1}^{\alpha}}\Lambda_{2}\Gamma^{\theta}$ & \\
\hline
$\bigoplus_{\alpha\in I_{2}}\Lambda_{1}\Gamma^{\alpha}\oplus
 \bigoplus_{\theta\in I_{1}}\Lambda_{1}\Gamma^{\theta}$
  & $\bigoplus_{\alpha\in I_{2}}\bigoplus_{\theta\in I_{1}^{\alpha}}\Lambda_{1}\Gamma^{\theta}$ & \\
\hline
$\Z^{L_{2}}\oplus \Z^{L_{1}} \oplus \Z^{L_{0}}$ & 
$\bigoplus_{\alpha\in I_{2}}\Z^{L_{1}^{\alpha}} \oplus
 \bigoplus_{\alpha\in I_{2}}\Z^{L_{0}^{\alpha}} \oplus 
 \bigoplus_{\theta\in I_{1}}\Z^{L_{0}^{\theta}} $ &
$\bigoplus_{\alpha\in I_{2}} \bigoplus_{\theta\in I_{1}^{\alpha}}\Z^{L_{0}^{\theta}}$ \\
\hline
\end{tabular}
}
\end{center}
\caption{First page of the Mayer-Vietoris spectral sequence for the homology of $\bA$.}
\label{tab:MVSS}
\end{table*}
Knowledge of $H_s(\bA)$ for $r=0$ and $1$ allows computation of the differentials, and the only torsion that can arise is that in the cokernel of the differential $d^1_{{1,2}}\co  E^1_{1,2}\to E^1_{0,2}$. This differential runs
$$
d^{1}_{{1,2}}\co  
  \underset{\alpha\in I_{2}}{\bigoplus} \left(\oplus_{\theta\in I_{1}^{\alpha}}
                \Lambda_{2} \Gamma^{\theta} \right)\rightarrow 
(\oplus_{\alpha\in I_{2}^{\ }} \Lambda_{2} \Gamma^{\alpha})
\,\oplus\,  
(\oplus_{\theta\in I_{1}^{\ }} \Lambda_{2} \Gamma^{\theta}) \notag
$$
and is described explicitly on each component $\Lambda_2\Gamma^\theta$ for $\theta\in I_1^\alpha$ via the canonical embeddings of the stabiliser subgroups $\Gamma^\theta$ into the corresponding $\Gamma^\alpha$ and $\Gamma^\theta$ in the target components. We obtain 

\begin{lemma} \label{TorinA}
For a rational projection pattern with $\nu=2$, the only torsion which arise in $H_*(\bA)$ is in $H_2(\bA)$ and is that which arises in the cokernel of the differential $d^{1}_{{1,2}}$.\qed
\end{lemma}

The cohomology calculations are formally dual to the above; note now that the important differential for torsion purposes runs $d_1^{{1,2}}\co   E_1^{0,2}\to E_1^{1,2}$, giving potential torsion in its cokernel, the group $E_2^{1,2}$, that is, in cohomological dimension 3. This corresponds, via the universal coefficient theorem, to the identification of the torsion in $H_2(\bA)$ with that in $H^3(\bA)$.

In all the icosahedral tilings we consider below, this torsion group is non-trivial. For larger values of $\nu$ the corresponding Mayer-Vietoris spectral sequences are similar, though have more rows. Again, torsion  can only arise from the cokernel of differentials $d^1_{1,s}$ where now $2\leqslant s\leqslant \nu$, and formulae for these differentials are given by the analogues of the description for $\nu=2$ above.

Finally we note the following useful observation concerning a criterion for the absence of torsion in the group $\coker(\alpha_*)$. We state it for general values of $\nu$.

\begin{lemma}\label{FullalphaZ}
For a codimension 3 rational projection pattern and $s>0$,  if $\coker \phi_s'$ is torsion free then $\Delta_s=0$ and hence $\coker(\alpha_*)$ is torsion free.
\end{lemma}

\smallskip\noindent{\bf Proof.}  
By Lemma \ref{FullalphaQ} we know that $\Delta_s=0$ when working over $\R$, and hence when working over $\Z$ the image of $\Delta_s$ can only be a torsion group. If $\coker \phi_s'$ is torsion free then $\Delta_s=0$ integrally. The result concerning $\coker(\alpha_*)$ now follows from  Lemma \ref{cokernels}\qed




\subsection{Codimension 3 examples with icosahedral symmetry}
\label{sec:ex3d}
We illustrate the above tools by considering the computations of
the cohomology of the four icosahedral tilings, the Danzer tiling
\cite{Dan89}, the Ammann-Kramer tiling \cite{AmmannKramer}, the
canonical $D_{6}$ tiling \cite{KP95} and the dual canonical $D_6$
tiling  \cite{KP95}, all of which are rational projection patterns.  
Preliminary results were announced in \cite{GHK}, but, as we note below, at least in the case of the Danzer tiling, the torsion component of the integral cohomology group $H^3(\MP)$ was incorrectly computed there, as possibly were also the corresponding computations of $H^3$ for the other three examples; the lower cohomology groups announced in \cite{GHK} are correct. In this section we give details of these computations. These examples also give a good overview of some of the different phenomena that can occur in the 
determination of torsion.

We start by describing the relevant lattices $\Gamma$ and families of
singular planes $\mathcal{W}$. In three dimensions, there are three 
inequivalent icosahedral lattices of minimal rank 6. 

The primitive
lattice $\Gamma_{P}$ is generated by a star of vectors pointing from
the center to the vertices of a regular icosahedron. We choose any
basis $e_{1},\ldots,e_{6}$ from this vector star. The lattice 
$\Gamma_{F}$ is then the sublattice of those integer linear combinations
of the $e_{i}$, whose coefficients add up to an even integer. 
The lattice $\Gamma_{I}$ is given by the $\Z$-span of the vectors in
$\Gamma_{P}$, and the additional vector $\frac12(e_{1}+\ldots+e_{6})$.
These lattices are analogues of the primitive, F-centered, and 
I-centered cubic lattices.\footnote{strictly speaking, $\frac12\Gamma_{F}$
is a centering of $\Gamma_{P}$} $\Gamma_{F}$ is an index-2 sublattice
of $\Gamma_{P}$, which in turn is an index-2 sublattice of $\Gamma_{I}$.
The action of the icosahedral group $A_{5}$ on the three lattices gives
rise to three integral representations, which are inequivalent under
conjugation in $GL_{6}(\Z)$. 

The singular planes of all four examples 
have special orientations, being perpendicular either to a 5-fold,
a 3-fold, or a 2-fold axis of the icosahedron (the latter are also parallel
to a mirror plane). Moreover, each $\Gamma$-orbit of singular planes 
contains a representative which passes through the origin. We therefore
define the families of planes $\mathcal{W}^{n}$, $n=5,3,2$, consisting
of all planes perpendicular to an $n$-fold axis, and passing through
the origin. The arrangements of singular planes of the icosahedral examples
are then given by the pair $(\Gamma_{P},\mathcal{W}^{2})$ for the 
Ammann-Kramer tiling, the pair $(\Gamma_{F},\mathcal{W}^{2})$ for 
the dual canonical $D_{6}$ tiling, the pair 
$(\Gamma_{F},\mathcal{W}^{5})$ for the the Danzer tiling, 
and the pair $(\Gamma_{F},\mathcal{W}^{5}\cup\mathcal{W}^{3})$ 
for the canonical $D_{6}$ tiling. Interestingly,
the sets $\Gamma_{P}+\mathcal{W}^{2}$ and $\Gamma_{F}+\mathcal{W}^{2}$
are invariant even under all translations from $\Gamma_{I}$, which
means that they are both equal to $\Gamma_{I}+\mathcal{W}^{2}$. In 
other words, the sets of singular planes of the Ammann-Kramer tiling 
and the dual canonical $D_{6}$ tiling are the same, only the lattices 
acting on it by translation are different. On the other hand, the sets
of singular planes of the Danzer tiling and the canonical $D_{6}$ tiling
have a lattice of translation symmetries which is equal to the lattice
$\Gamma_{F}$ they are constructed from. With these data is it now 
straightforward to evaluate the formul\ae\ of
Theorem~\ref{RatFormulae} 
for the ranks of the rational homology groups. 
The results are summarized in Table~\ref{tab:ex}. As can be seen, compared 
to previously published results the rational ranks of $H_{0}$ and $H_{1}$ of 
the Ammann-Kramer tiling and the dual canonical $D_{6}$ tiling have 
been increased by 1, in agreement with Kalugin \cite{Kal}, whereas 
all other rational ranks remain the same.

Next, we discuss the determination of torsion, which is potentially
non-trivial only for $H^2$ and $H^3$, as noted by Proposition \ref{notor}. The torsion in $H^2$ is relatively straightforward, and can be computed, as in Theorem \ref{codim3}, as the torsion in the cokernel of 
$$\phi_2\co  \oplus_{\alpha\in I_2}\Lambda_4\Gamma^\alpha\to\Lambda_4\Gamma\,.$$    
Both the Ammann-Kramer and dual canonical $D_6$ tilings have 2-torsion in $H^2(\MP)$ arising from this map; $\coker\phi_2'$ is torsion free for the Danzer and canonical $D_6$ tilings, making their $H^2$ groups free abelian.

This situation concerning the groups $H^3$ is as follows. In each of the first three tilings $\coker\phi_1'$ is torsion free and so by Lemma \ref{FullalphaZ} the map $\Delta_1$ is zero. This is not so for the dual canonical $D_6$ tiling where additional geometric computation (such as that described in \cite{GHK}) is needed to deduce that nevertheless $\coker(\alpha_*)$ is still free. By Lemma \ref{cokernels}, this implies that the map $\Delta_1$ is non-trivial, having image the whole of the torsion subgroup 
$\Z_2^6$ of $\coker\phi_1'$.  

In every case there are 2-torsion components in $\coker\phi_1''$.  For the first three examples there is no other torsion arising in the computation via Diagram (\ref{diagramn=3}), but there remains an extension problem of the form
$$0\to\mbox{ free abelian group }\to \coker\phi_1\to\coker\phi_1''\to0$$
which has more than one potential solution.

For the dual canonical $D_6$ tiling, the situation is more complicated. There is torsion, $\Z_2^7$, in $\coker\phi_1''$, and an extension problem to decide its lift to $\coker\phi_1$, but there is also torsion, $\Z_2^{15}$,  in $\ker\phi_0$, arising as the torsion in $\ker\phi_0'$. Then $H^3$ is given by a further extension (\ref{extn}) running
\begin{equation}\label{DCD6}0\to\coker\phi_1\to H^3(\Omega)\to\Z^{328}\oplus\Z_2^{15}\to0\,.\end{equation}
There is no direct way of solving such extension problems without some additional geometric input, for example via the geometric realisation of the rational projection pattern and computation of the corresponding torus arrangement $\bA$. Using the Mayer-Vietoris computation of $H_*(\bA)$, it can be shown that, for this example, at least one of these extension is non-trivial.

Clearly each example at this point must be handled on a case by case basis. The Danzer tiling, arguably the simplest of these four, involves sufficiently small cell complexes that a complete solution is available. Computations using either Diagram (\ref{diagramn=3}) or the exact sequence (\ref{Ples}) yield an extension problem for $H^3$ with a single $\Z_2$ in the quotient: there is a simple dichotomy, that $H^3$ is either $\Z^{20}\oplus\Z_2$ (the trivial extension), or is $\Z^{20}$ (the non-trivial one). However, for this tiling, a modified version of the Anderson-Putnam complex gives $H^3(\MP)$ as the direct limit of the cohomology $H^3(K)$ of a certain cell complex under an iterated self map $f\colon K\to K$. Machine computation shows that there is no torsion in $H^3(K)$, and consequently there can be none in the direct limit $H^3(\MP)=\lim_\to H^3(K)$. This determines the extension problem: it is the non-trivial one.

\begin{cor} The integral cohomology of the Danzer tiling is given by
$$H^0(\MP)=\Z\qquad H^1(\MP)=\Z^7\qquad H^2(\MP)=\Z^{16}\qquad H^3(\MP)=\Z^{20}\,.$$
\end{cor}\qed




\begin{table*}
\begin{center}
{
\renewcommand{\arraystretch}{1.4}
\begin{tabular}{|c|c|c|c|c|c|c|l|}
\hline
&$H^3\otimes\Q$                          & $H^2$                      & $H^1$     & $H^0$ & $t_1'$ & $t_1''$   & $t_0'$  \\
\hline    Danzer     &
$\Q^{  20}$             & $\Z^{ 16}$                 & $\Z^{ 7}$ & $\Z$  &     0 & $\Z_2$  & 0       \\
\hline    Ammann-Kramer &
$\Q^{ 181}$             & $\Z^{ 72}\oplus\Z_2$        & $\Z^{12}$ & $\Z$  &  0 & $\Z_2$  & 0        \\
\hline     canonical $D_6$     &
$\Q^{ 205}$         & $\Z^{ 72}$                 & $\Z^{ 7}$ & $\Z$  &    0 & $\Z_2^2$& 0     \\
\hline      dual canonical $D_6$    &
$\Q^{ 331}$ & $\Z^{102}\oplus\Z_2^4\oplus\Z_4$ & $\Z^{12}$ & $\Z$  &    $\Z_2^6$& $\Z_2^7$   & $\Z_2^{15} $    \\
\hline
\end{tabular}
}\smallskip
\end{center}
\caption{\cite{GHK} Integral cohomology $H^2$, $H^1$ and $H^0$, and rational $H^3$ of icosahedral tilings from the literature. Also indicated are details of the torsion arising at various points in the calculation via diagram (\ref{diagramn=3}). We use the notation that $t_1'$, $t_1''$ and $t_0'$ denote the torsion components of $\coker\phi_1'$, $\coker\phi_1''$ and $\ker\phi_0'$ respectively.}
\label{tab:ex}
\end{table*}

\subsection{$K$-theory and general codimension}\label{Gen}
We conclude with some remarks on the general codimension case and the consequences of our work for the $K$-theory of cut \& project patterns. The first observation is an extension of the result of Corollary \ref{PThigh}.

\begin{prop} For a general dimension $d$, codimension $n$, rational projection tiling, the groups $H^s(\Omega)$ are free abelian of rank $N\choose s$ if $s<\nu-1$, but the possibility of torsion exists for all $s\geqslant\nu-1$. For these values of $s$ the torsion subgroup contains the torsion part of
$$\coker\left\{\bigoplus_\alpha\Lambda_{N-s}\Gamma^\alpha
\buildrel\alpha_*\over\longrightarrow \Lambda_{N-s}\Gamma\right\}$$
and if $s<2\nu-1$ then this is precisely the torsion term.
\end{prop}

\smallskip\noindent{\bf Sketch of Proof.} By Corollary \ref{Ples} we have an exact sequence
$$\cdots \to H_{N-s}(\bA)\buildrel\alpha_*\over\longrightarrow H_{N-s}(\bT)
\longrightarrow H^s(\MP)\longrightarrow H_{N-s-1}(\bA)\to\cdots\,.$$
If $s<\nu-1$, then $N-s-1>(n-1)\nu$ and $H_{N-s-1}(\bA)=0$ since $\bA$ is constructed as the union of $(n-1)\nu$-dimensional tori. This proves the first part of the statement, a restatement of Corollary \ref{PThigh}. The exact sequence also shows that the cokernel term in the statement clearly injects in $H^s(\MP)$, which forms the second part, and the final observation follows from a more detailed computation of $H_*(\bA)$. In brief, in the range $s<2\nu-1$, the groups $H_{N-s-1}(\bA)$ are free abelian; this follows from a Mayer-Vietoris spectral sequence computation of $H_*(\bA)$ as in \cite{Kal}, and uses the observation that the highest homological dimension of intersection of the component $(n-1)\nu$-tori making up $\bA$ is of dimension $(n-2)\nu$.\qed

\smallskip
Examples lead us to conjecture that where this proposition indicates that there might be torsion, examples can be found where there is torsion.

\smallskip
We turn to consider the various forms of $K$-theory used in the study of aperiodic tilings. Initial interest was in the $K$-theory of various noncommutative $C^*$-algebras associated to the tilings \cite{Bellissard}. At the level of graded abelian groups, the values of this $K$-theory is (modulo regrading) the same as the topological $K$-theory of the tiling space, which we denote $K^*(\MP)$ \cite{FoHu}. The Atiyah-Hirzebruch spectral sequence (AHSS) provides a way of computing the $K$-theory of $\MP$ from its cohomology, and in the absence of torsion in $H^*(\MP)$ it is a standard fact that this spectral sequence collapses and the $K$-theory is just the direct sum of the cohomology groups, $K^0(\MP)=\oplus_sH^{2s}(\MP)$ and $K^1(\MP)=\oplus_sH^{2s-1}(\MP)$. This is not necessarily so in the presence of torsion. However, for small values of $d$ the situation remains straightforward. Standard topological arguments with the AHSS, using characterisation of the possible differentials yield

\begin{prop}
For a dimension $d\leqslant 3$ rational projection pattern, $K^0(\MP)=\oplus_sH^{2s}(\MP)$ and $K^1(\MP)=\oplus_sH^{2s-1}(\MP)$.
\end{prop}

\smallskip\noindent{\bf Sketch of Proof.} The smallest possible non-zero differential that can lead to the failure of this result is $d_3\co  H^s(\MP)\to H^{s+3}(\MP)$ and there must be 2-torsion in $H^{s+3}(\MP)$. However, as $H^0(\MP)=\Z$ and represents the connectivity element, $d_3$ can be non-trivial only for positive values of $s$; as we need a non-trival value of $H^{s+3}(\MP)$ for $d_3\not=0$, we need the cohomology of $\MP$ to be non-zero in some dimensions at least 4. This cannot happen if the tiling is only 3 dimensional. 

As we know that $\nu$ must be an integer, the codimension of the pattern is 3 or less, and torsion can only exist in cohomological dimensions 2 or 3; this is enough to be sure that there are no extension problems and the $K$-theory splits as the direct sum as claimed.\qed

The possibility nevertheless exists for tilings in higher dimensions demonstrating significant differences between their cohomology and $K$-theory.

\section{Appendix: realisation of $H^*(\MP)$ as group (co)homology}
We  sketch briefly here the argument that $H^*(\MP)$ can be realised as the group cohomology $H^*(\Gamma;C_n)$ as noted at the start of Section \ref{HomAlg}. This identification lies at the heart  of the approach of \cite{FHKcmp,FHKmem} and we direct the reader to \cite{FHKmem} in particular for full details. It also enables us to justify the identification of the homomorphisms $\mu^*\co  H^*(\bT)\to H^*(\MP)$ and $m_*\co  H^*(\Gamma; \uT_*)\to H^*(\Gamma;C_n[n])$ as mentioned in Remark \ref{identMU}.

We recall from Section \ref{Sect2} that $\MP$ can be considered as the completion of $q(NS)\subset\bT$, the image of the non-singular points in the $N$-torus, with respect to the \met.

The translation action of $E^\|$ on $\Eup$ passes to an action on $\bT=\Eup/\Gammaup$. Moreover, this action preserves $S$ and $N\! S$ and the resulting action on $q(N\! S)$ is continuous in the \met. Thus the dynamical system $(q(N\! S),\Epar)$ extends to the completion yielding the dynamical system $(\MP,E^\|)$. Since $(q(N\! S),E^\|)$ is a sub-system of a flow on the torus $\bT$, the system 
$(\MP,E^\|)$ is the flow on the mapping torus of a $\Z^d$ action. 

To make this explicit, choose a splitting $\Gamma=\Gamma_1\oplus\Gamma_2$ where $\Gamma_1$ is of rank $n$, and $\Gamma_2$ is of rank $d$. Then $\Gamma_1$ spans a linear space $F\subset\Eup$ and we can pick a fundamental domain $X\subset F$. Let $X_c$ be the completion of $q(X\cap N\! S)$ in the \met\ and $F_c$ the corresponding completion of $q(F\cap N\! S)$. Then $\MP$ is the mapping torus of the action of $\Gamma_2\cong\Z^d$ on $X_c$ given on $q(x)\in q(X\cap N\! S)$ by $\gamma\cdot q(x) = q(x-\pi_F(\gamma))$, where $\pi_F$ is the projection onto $F$ along $\Epar$. 

Identifying $\MP$ with this mapping torus, it follows that the Cech cohomology of $\MP$ can be identified with $H^*(\Gamma_2;C(X_c,\Z))$ the group cohomology of $\Gamma_2$ with values in the representation module $C(X_c,\Z)$, the continuous $\Z$-valued functions on $X_c$, defined by this action. This can be seen, for example, by viewing the mapping torus as a fibre bundle over the $d$-torus $B\Gamma_2$ with fibre $X_c$: then the Serre spectral sequence of the bundle has $E_2=E_\infty$-term $H^*(\Gamma_2;H^*(X_c))$, and the Cech cohomology $H^*(X_c)$ of the totally disconnected space $X_c$ is precisely $C(X_c,\Z)$.  Moreover, as $\Gamma_1$ acts freely on $C(F_c,\Z)$ by the analogous action, we can identify $H^*(\Gamma_2;C(X_c,\Z))$ with $H^*(\Gamma;C(F_c,\Z))$.
Finally, the assumption of total irrational position allows us to identify $C(F_c,\Z)$ with $C(\Eperp_c,\Z))$, as $\Gamma$-modules, where $\Eperp_c$ is the completion of $\Eperp\cap N\!S$ in the \met, as in Section \ref{HomAlg}. Summing up, we have the string of identifications
$$H^*(\MP)\cong H^*(\Gamma_2;C(X_c,\Z))\cong H^*(\Gamma;C(F_c,\Z))\cong H^*(\Gamma;C(\Eperp_c,\Z))\,.$$
As noted after Definition \ref{C-tope}, we may identify $C(\Eperp_c,\Z))$ with the $\Gamma$-module $C_n$. The equivalence $H^s(\Gamma;C_n)\cong H_{d-s}(\Gamma;C_n)$ in (\ref{FundEquiv}) follows from the Poincar\'e duality property of the group $\Gamma_2$.

We turn to the map $\mu\co \MP\to\bT$. The identifications above allow us to view $\MP$ cohomologically as the fibre bundle over $B\Gamma$ with fibre $\Eperp_c$; in the same way, we can replace the space $\bT$ by the bundle  over $B\Gamma$ with fibre $\Eperp$ and the homomorphism $\mu^*$ is represented in cohomology by $H^*(\Gamma;C(\Eperp;\Z))\to H^*(\Gamma;C(\Eperp_c;\Z))$ induced by the natural quotient $\Eperp_c\to\Eperp$. By the analogous argument along the lines of the proof of Theorem \ref{identification}, this is the same map in cohomology as induced in $H^*(\Gamma; -)$ by the map $\uT_*\to C_n[n]$.

\end{document}